\documentclass[10pt]{article}

\usepackage[margin=1.0in]{geometry}
\usepackage{amsthm, amsfonts,amssymb,euscript}
\usepackage{amsmath}
\usepackage{enumerate}
\usepackage{bm}
\usepackage{amssymb}
\usepackage{amsthm}
\usepackage{graphicx}
\usepackage{caption}
\usepackage{subcaption}
\usepackage{amsfonts}
\usepackage{float}
\usepackage{color}
\usepackage{slashed}
\usepackage[affil-it]{authblk}
\usepackage{mathrsfs}
\usepackage{upgreek}
\usepackage{array,makecell}
\usepackage{boldline,multirow}

\usepackage[usenames,dvipsnames,svgnames,table]{xcolor}
\usepackage[colorlinks=true]{hyperref}
\hypersetup{linkcolor=BrickRed, urlcolor=green, citecolor=blue, linktoc=page}

\numberwithin{equation}{section}

\newtheorem*{proposition*}{Proposition}
\newtheorem*{theorem*}{Theorem}
\newtheorem*{conjecture*}{Conjecture}
\newtheorem*{claim*}{Claim}
\newtheorem*{lemma*}{Lemma}
\newtheorem*{corollary*}{Corollary}
\newtheorem{theorem}{Theorem}[section]

\newtheorem*{definition*}{Definition}
\newtheorem{definition}{Definition}[section]
\newtheorem*{assumption*}{\mathcal{A}ssumption}

\newtheorem*{remark*}{Remark}
\newtheorem{remark}{Remark}[section]

\allowdisplaybreaks

\begin{document}

\title{Supersonic flow past an airfoil}

\date{}
\author[1]{Yannis Angelopoulos}
\affil[1]{\small The Division of Physics, Mathematics and Astronomy, Caltech,
1200 E California Blvd, Pasadena CA 91125, USA}

\normalsize

\maketitle
\begin{abstract}
We consider a supersonic flow past an airfoil in the context of the steady, isentropic and irrotational compressible Euler equations. We show that for appropriate data, either a shock forms, in which case certain derivatives of the solution blow up, or a sonic line forms, which translates to the fact that the hyperbolic character of the equations degenerates. 
\end{abstract}

\section{Introduction}
In this paper we make several remarks on the study of shock formation in the case of supersonic flows around airfoils. This is a problem coming from physics and engineering concerning the properties of supersonic flight.

\begin{figure}[H]
\begin{center}
\includegraphics[width=6cm]{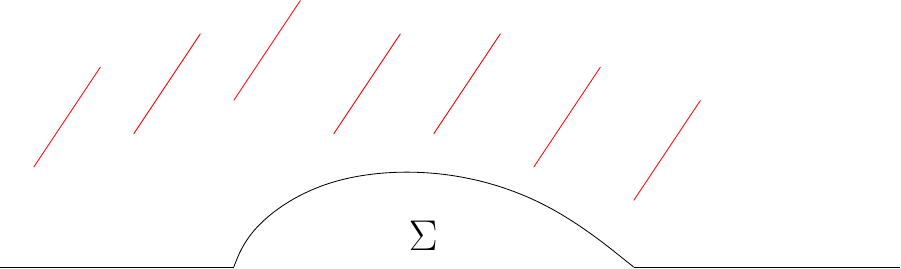}
\vspace{0.4cm}
\caption{\label{supersonic-figure1}The initial hypersurface $\Sigma$, and the dotted area in red where the solution is constructed, and a shock and a sonic line form.}
\end{center}
\end{figure}

Our setup is that of the 2-dimensional steady, isentropic and irrotational Euler equations as an initial value problem, where the initial data are given on a hypersurface, part of which models an airfoil. On the initial hypersurface we assume that the Mach number is bigger than one, making our equations hyperbolic. For data of appropriate size we show that either a shock or a sonic line forms. We are also able to have a precise estimate for the time where the shock forms and where the sonic line forms relative to the size of the initial data. A rough version of the main result of the paper is the following:

\begin{theorem}\label{thm:mainrough}
Consider a solution of the steady, isentropic and irrotational Euler equations in two dimensions (denoted by variables $x$ and $y$) given in the form of a potential flow $\phi$ emanating from data set on an initial smooth hypersurface $\Sigma$ as in figure \ref{supersonic-figure1}, where we assume that
$$ \sqrt{ ( \partial_y \phi )^2 + ( \partial_x \phi )^2 } > q_{cr} , $$
a condition that makes the second order quasilinear equation satisfied by $\phi$ initially hyperbolic.

Then in the upper half space relative to $\Sigma$ we construct the maximal domain of existence and uniqueness of the potential flow $\phi$, and on the boundary of that domain either a shock forms, i.e. as we approach the boundary we have that $\partial^2 \phi \rightarrow \infty$, or a sonic line forms, i.e. as we approach the boundary we have that 
$$\sqrt{ ( \partial_y \phi )^2 + ( \partial_x \phi )^2 } \rightarrow q_{cr}. $$
\end{theorem}

\vspace{-0.8cm}

\begin{remark}
Let us note that the boundary of the maximal domain of existence and uniqueness of the solution constructed here can contain points where a shock forms or points where a sonic line forms. The full maximal boundary is expected to contain both, such a boundary is sometimes referred to in the literature as a \textbf{shock polar}. 
\end{remark}

As far as we know, the problem of the two dimensional steady, isentropic and irrotational flow has been studied extensively from a mathematical viewpoint, but only in the subsonic and transonic cases. For the subsonic case (which involves the study of a quasilinear elliptic equation) see the work of Bers \cite{bers-subsonic1}, \cite{bers-subsonic2}, \cite{bers-book} and Shiffman \cite{shiffman}, while for the transonic case (where the equation changes character from hyperbolic to elliptic) we refer to the important work of Morawetz \cite{morawetz1}, \cite{morawetz2}, \cite{morawetz3} (see also the review article \cite{morawetz-ams} and the companion paper of the author \cite{transonic1})). 

For the supersonic case, there have been numerical results and heuristical arguments that point towards the formation of shocks and sonic lines, but no rigorous mathematical proof. For a reference on numerical and heuristical results one can consult the book of Courant and Friedrichs \cite{courant-friedrichs} (see also the lecture notes of Morawetz \cite{morawetz-notes}).

On the other hand, during the last 15 years, there has been significant progress on the shock formation problem for hyperbolic equations. The main reference is the monumental work of Christodoulou \cite{DC07} where he constructs shocks for three dimensional relativistic fluids. For the non-relativistic case see \cite{christ-miao}. Alinhac has also obtained weaker results on shock formation for quasilinear wave equations, see \cite{alinhac-blow1}, \cite{alinhac-blow2} (for even earlier results see the pioneering work of John \cite{john-book}). The methods of Christodoulou have been used to show shock formation for quasilinear wave equations, see the book of Speck \cite{speck-book}. Since then the same methods have been used to show shock formation in several different scenarios, always for the time-dependent Euler equations, see \cite{jaredgustavjonwillie}, \cite{hksw}, \cite{jaredjonathan} to mention a few. Speck has also used Christodoulou's method to show degeneration of hyperbolicity for certain quasilinear wave equations \cite{jared-hyperbolic}. 

In this work we also use the methods of Christodoulou. A natural question is in the case where a shoci forms, whether the solution can be continued beyond the shock and reach the sonic line, hence constructing a shock polar. The problem of shock development has been considered in the important works of Christodoulou-Lisibach \cite{christ-lisibach} and Christodoulou \cite{christ-shockdev} (see also the work of Majda \cite{majda-1}, \cite{majda-2}), but the main issue in the present situation is the formation of the sonic line. Note that we cannot expect to continue the solution past the sonic line as such a problem is in general ill-posed (see for instance \cite{metivier-nonlinear}, \cite{lerner-morimoto-xu}, \cite{lerner-nguyen-texier}, \cite{transonic1}).

\section{The steady, isentropic and irrotational Euler equations in 2D}

We want to study the 2-dimensional irrotational, isentropic and \textit{steady} Euler flow past a symmetric airfoil $\mathcal{P}$ (that will be defined later). 

Consider first an isentropic and irrotational flow on $\mathbb{R} \times \mathbb{R}^2$ for which we have the following sets of equations:
\begin{equation}\label{eq:euler}
\begin{cases}
\partial_t \rho + \partial_{x_1} ( \rho u_1 ) + \partial_{x_2} ( \rho u_2 ) = 0 , \\
\partial_t u_i + u_j \cdot \nabla_{x_1 ,x_2 } u_i - \frac{1}{\rho} \partial_{x_i} p = 0 \mbox{  for $i, j \in \{1,2\}$},
\end{cases}
\end{equation}
where $\rho$ is the density, $(u_1 , u_2 )$ is the velocity vector, and $p \doteq p (\rho )$ is the pressure. Let us denote by $q$ the speed, i.e.
$$ q = \sqrt{u_1^2 + u_2^2} , $$
and by $c$ the speed of sound given by
$$ c^2 = - \frac{\rho q}{\rho' (q)} , $$
something that follows from Bernoulli's law.

Now we make the assumption that the flow is steady, i.e. independent of time. We note also that as $p$ depends only on $\rho$ it can be eliminated from \eqref{eq:euler}, and by the assumption that the flow is irrotational as we have that
$$\partial_{x_2} u_1 - \partial_{x_1} u_2 = 0 , $$
we can find a function $\phi : \mathbb{R}^2 \rightarrow \mathbb{R}$ such that
$$ u_1 = \partial_{x_1} \phi , \quad u_2 = \partial_{x_2} \phi . $$

Setting $x_1 = x $ and $x_2 = y$, the potential function $\phi$ then satisfies the following equation:
\begin{equation}\label{eq:phi}
\left( 1 - \frac{( \partial_{x} \phi )^2}{c^2} \right) \partial_{xx}^2 \phi - \frac{2 \partial_x \phi \partial_y \phi}{c^2} \partial_{xy}^2 \phi + \left( 1 - \frac{ ( \partial_y \phi )^2}{c^2} \right) \partial_{yy}^2 \phi = 0 .
\end{equation}
Letting 
$$ M = \frac{q}{c} = \frac{\sqrt{(\partial_x \phi)^2 + (\partial_y \phi)^2}}{c} , $$
we note that equation \eqref{eq:phi} is \textbf{elliptic} if $M < 1$ and \textbf{hyperbolic} if $M >1$. We call the flow \textit{subsonic} if equation \eqref{eq:phi} is elliptic, and we call it \textit{supersonic} if equation \eqref{eq:phi} is hyperbolic.

\begin{remark}
Another way to study the problem of the irrotational, isentropic and steady flow in $1+2$ dimensions (or just two dimensions as time is eliminated) is to use the equation
$$ \partial_x ( \rho u_1 ) + \partial_y ( \rho u_2 ) = 0 , $$
in order to introduce a stream function $\psi$ such that
$$ \rho u_1 = \partial_y \psi , \quad \rho u_2 = - \partial_x \psi . $$

For $\psi$ we can use then the so called \textit{hodograph transformation}. Let 
$$ \vartheta = \arctan \left( \frac{u_2}{u_1} \right) , $$
and let us consider the change of coordinates
$$ (x,y) \rightarrow ( \vartheta, q ) . $$
This is a dynamic change of variables, for which we note that
$$ d\phi = u_1 dx + u_2 dy , \quad \frac{1}{\rho} d\psi = -u_2 dx + u_1 dy , $$
which implies that for $z = x + iy$ we have that
$$ d\phi + i \frac{1}{\rho} d\psi = q e^{-i\vartheta } dz \Rightarrow dz = \frac{e^{i\vartheta}}{q} \left( d\phi + i \frac{1}{\rho} d\psi \right) . $$
We can now consider the partial derivatives of $z$ and we have that
$$ \frac{\partial z}{\partial \vartheta} = \frac{e^{i\vartheta}}{q} \left( \frac{\partial \phi}{\partial q} + i \frac{1}{\rho} \frac{\partial \psi}{\partial q} \right) , \quad \frac{\partial z}{\partial \vartheta} = \frac{e^{i\vartheta}}{q} \left( \frac{\partial \phi}{\partial \vartheta} + i \frac{1}{\rho} \frac{\partial \psi}{\partial \vartheta} \right) .$$
Since $\frac{\partial^2 z}{\partial q \partial \vartheta} =  \frac{\partial^2 z}{\partial \vartheta \partial q}$, we have the following \textit{hodograph equations}:
$$ \frac{\partial \phi}{\partial q} = q \frac{\partial \psi}{\partial \vartheta} \frac{\partial}{\partial q} \left( \frac{1}{q \rho} \right) , \quad \frac{\partial \phi}{\partial \vartheta}= \frac{q}{\rho} \frac{\partial \psi}{\partial q} . $$
Using Bernoulli's Theorem (that $\frac{dp}{dq} = - q \rho$) and that $\frac{\partial^2 \phi}{\partial q \partial \vartheta} = \frac{\partial^2 \phi}{\partial \vartheta \partial q}$ we have the following second order equation for $\psi$:
\begin{equation}\label{eq:psi_h}
( 1 - M^2 ) \frac{\partial^2 \psi}{\partial \vartheta^2} + q^2 \frac{\partial^2 \psi}{\partial q^2} + q ( 1+M^2 ) \frac{\partial \psi}{\partial q} = 0 ,
\end{equation}
for $M = \frac{q}{c}$ as before.
\end{remark}

\section{The initial boundary problem for airfoils and some geometric assumptions}\label{data}

We will study equation \eqref{eq:phi} with data on 
$$\Sigma \doteq \mathcal{P} \cup \{ (x,y) \quad | \quad y=0 , |x| \geq 1 \} , $$
where $\mathcal{P}$ is an airfoil that is defined as the graph of a smooth concave function $f$ such that 
$$ f(x) > 0 \mbox{  for $|x| < 1$,  } f(1) = f(-1) = 0 . $$
Note that we have that
$$ \mathcal{P} \doteq \{ (x,y) \quad | \quad |x| < 1 , y = f(x) \} . $$

We will examine solutions of \eqref{eq:phi} only in the domain
$$ \mathcal{D} \doteq \{ (x, y ) \quad | \quad y \geq 0 \mbox{  for $|x| \geq 1$, and  } y \geq f(x) \mbox{  for $|x| < 1$} \}. $$

On $\Sigma$ we let
$$ \left. \phi \right|_{\Sigma} \in C^{\infty} (\mathbb{R} ) , \quad \left. \frac{\partial \phi}{\partial n} \right|_{\Sigma} = 0 . $$
We also let
$$ \lim_{|x| \rightarrow \infty} q(x) = K \mbox{  such that  } \lim_{|x| \rightarrow \infty} M(x) > 1 . $$
Such a problem can be seen as an initial-boundary value problem, where we assign a supersonic speed at infinity.

As noted before, if in $\mathcal{D}$ we have that $M > 1$ then equation \eqref{eq:phi} is a quasilinear hyperbolic equation in $\mathcal{D}$ in $1+1$ dimensions. If the condition $M >1$ is satisfied by the initial data then we can treat \eqref{eq:phi} as a hyperbolic equation, but this might not be the case.

For a given density function we note that there exists a number $q_{cr}$ (the \textit{critical speed}) such that
$$ M > 1 \mbox{  if $q > q_{cr}$, and  } M <1 \mbox{  if $q < q_{cr}$,} $$
i.e. equation \eqref{eq:phi} is hyperbolic for $q > q_{cr}$ and elliptic for $q < q_{cr}$.

\section{The purely hyperbolic regime}\label{hyperbolic}

\subsection{The equations}

Since we will work with data that will guarantee that $M > 1$ in the domain where we will show existence of the solution, equation \eqref{eq:phi} can then be treated as a quasilinear hyperbolic equation of the form:
$$ \bar{g}^{\alpha \beta} ( \partial \phi ) \partial^2_{\alpha \beta} \phi = 0 \mbox{  for $\alpha , \beta \in \{x,y\}$} , $$
where
$$ \bar{g}^{yy} = 1 - \frac{( \partial_y \phi )^2}{c^2} , \quad \bar{g}^{xy} = - \frac{( \partial_x \phi ) ( \partial_y \phi )}{c^2} , \quad \bar{g}^{xx} = 1 - \frac{( \partial_x \phi )^2}{c^2} .$$

Rewriting equation \eqref{eq:phi} as 
$$ ( c^2 - ( \partial_x \phi )^2 ) \partial_{xx}^2 \phi - 2 ( \partial _x \phi ) ( \partial_y \phi ) \partial_{xy}^2 \phi + ( c^2 - ( \partial_y \phi )^2 ) \partial_{yy}^2 \phi = 0 , $$
we multiply the equation above with $- \frac{1}{c^2 - ( \partial_y \phi )^2}$ (as suggested in \cite{speck-book}), and after differentiating it by $\partial_y$ and $\partial_x$ for
$$ \Phi_1 \doteq \partial_y \phi , \quad \Phi_2 \doteq \partial_x \phi , $$
we get the following system for $\Phi = ( \Phi_1 , \Phi_2 )$:
\begin{equation}\label{eq:Phi}
\begin{cases}
\Box_{g( \Phi )} \Phi_1 = \mathcal{N}_1 ( \Phi_1 , \Phi_2 ) , \\
\Box_{g (\Phi )} \Phi_2 = \mathcal{N}_2 ( \Phi_1 , \Phi_2 ) ,
\end{cases}
\end{equation}
where
\begin{equation*}
g^{-1} ( \Phi ) =
\begin{bmatrix}
g^{00} & g^{01} \\
g^{10} & g^{11}
\end{bmatrix} 
=
\begin{bmatrix}
-1 & \frac{ \Phi_1 \cdot \Phi_2}{c^2 - \Phi_1^2} \\
\frac{ \Phi_1 \cdot \Phi_2}{c^2 - \Phi_1^2} & - \frac{ c^2- \Phi_2^2}{c^2 - \Phi_1^2} 
\end{bmatrix} ,
\end{equation*}
which implies that
\begin{equation*}
g ( \Phi ) =
\begin{bmatrix}
g_{00} & g_{01} \\
g_{10} & g_{11}
\end{bmatrix} 
=
\begin{bmatrix}
- \frac{(c^2 - \Phi_1^2 )(c^2 - \Phi_2^2 )}{c^2 ( c^2 - \Phi_1^2 - \Phi_2^2 )} & \frac{(c^2 - \Phi_1^2) \Phi_1 \cdot \Phi_2}{c^2 ( c^2 - \Phi_1^2 - \Phi_2^2 ) } \\
\frac{ ( c^2 - \Phi_1^2 ) \Phi_1 \cdot \Phi_2}{c^2 ( c^2 - \Phi_1^2 - \Phi_2^2 )} & - \frac{ ( c^2- \Phi_1^2 )^2}{c^2 ( c^2 - \Phi_1^2 - \Phi_2^2 )} 
\end{bmatrix} ,
\end{equation*}
and where $\mathcal{N}_1$, $\mathcal{N}_2$ are the following nonlinear terms:
\begin{align*}
\mathcal{N}_1 ( \Phi_1 , \Phi_2) & = \frac{\partial}{\partial x} \left( \frac{\Phi_1 \cdot \Phi_2}{c^2 - \Phi_1^2} \right) \partial_y \Phi_1 - \frac{\partial}{\partial y} \left( \frac{\Phi_1 \cdot \Phi_2}{c^2 - \Phi_1^2} \right) \partial_x \Phi_1 \\ & + \frac{\partial}{\partial y} \left( \frac{c^2 - \Phi_2^2}{c^2 - \Phi_1^2} \right) \partial_x \Phi_2 - \frac{\partial}{\partial x} \left( \frac{c^2 - \Phi_2^2}{c^2 - \Phi_1^2} \right) \partial_x \Phi_1 \\ &+ \frac{1}{\sqrt{-\det (g)}} [ \partial_{\alpha}  ( \sqrt{-\det (g)} ) ] \cdot  g^{\alpha \beta} \partial_{\beta} \Phi_1 ,
\end{align*}
\begin{align*}
\mathcal{N}_2 ( \Phi_1 , \Phi_2 ) & = \frac{\partial}{\partial y} \left( \frac{\Phi_1 \cdot \Phi_2}{ c^2 - \Phi_1^2} \right) \partial_x \Phi_2 - \frac{\partial}{\partial x} \left( \frac{\Phi_1 \cdot \Phi_2}{ c^2 - \Phi_1^2} \right) \partial_y \Phi_2 \\ & + \frac{1}{\sqrt{-\det (g)}} [ \partial_{\alpha}  ( \sqrt{-\det (g)} ) ] \cdot  g^{\alpha \beta} \partial_{\beta} \Phi_2 .
\end{align*}

\subsection{Coordinate systems}
Let us define the eikonal function $u$ as the solution of the equation:
\begin{equation}\label{eikonal}
\begin{split}
& g^{\alpha \beta} \partial_{\alpha} u \partial_{\beta} u = 0 , \quad \partial_y u > 0 ,
\\ & \left. u \right|_{\Sigma} = 1-x \mbox{  for $|x| > 1$ and extended smoothly in the remaining part of $\Sigma$.}
\end{split}
\end{equation}
Define the vector field $L$ with components:
\begin{equation}\label{def:lmu}
L^{\alpha} = - \mu g^{\alpha \beta} \partial_{\beta} u , 
\end{equation}
where
\begin{equation}\label{def:mu}
\mu = - \frac{1}{g^{0\alpha} \partial_{\alpha} u} ,
\end{equation}
is the \textit{inverse foliation density}, which is the key quantity controlling the shock formation phenomenon.

It can be checked that $L$ is null, i.e.
$$ g(L , L ) = 0 . $$
Define also the vector field $\underline{L}$ as one satisfying
$$ g( L , \underline{L} ) = -2 , $$
and also the vector field
$$ X = \mu \frac{\underline{L} - L}{2} , $$
with components
$$ X^{\alpha} = - L^{\alpha} - g^{0\alpha} \mbox{  for $\alpha \in \{0,1\}$.} $$
Note that
$$ X^0 = 0 , \quad X^1 = - L^1 + \frac{\Phi_1 \Phi_2}{c^2 - \Phi_1^2} . $$
We also define the auxiliary vector field
$$ \check{X} = \frac{\underline{L} - L}{2} = \mu^{-1} X . $$
The Christoffel symbols can be expressed as
$$ \Gamma_{\alpha \gamma \beta } = \frac{1}{2} ( \partial_{\alpha} g_{\beta \gamma} + \partial_{\beta} g_{\alpha \gamma} - \partial_{\gamma} g_{\alpha \beta} ) \doteq \frac{1}{2} ( H_{\beta \gamma}^{\lambda} \partial_{\alpha} \Phi_{\lambda} + H_{\alpha \gamma}^{\lambda} \partial_{\beta} \Phi_{\lambda} - H_{\alpha \beta}^{\lambda} \partial_{\gamma} \Phi_{\lambda} ) , $$
for $\alpha , \beta , \gamma \in \{0,1\}$, where
$$ \partial_0 = \partial_y , \quad \partial_1 = \partial_x , $$
and where the functions $H_{\alpha \beta}^{\lambda}$, $\alpha , \beta \in \{0,1\}$, $\lambda \{1,2\}$ are defined from the aforementioned inequality. 
Finally, for $y = t$, and $u$ as defined before, we consider the coordinate system $(t,u)$ that we will call \textit{dynamical coordinates} due to their dependence on the solution itself. For the change of variables map 
$$ \mathcal{J} ( t,u) = (y,x) , $$
we have the following computation for the Jacobian:
\begin{equation}\label{jacobian}
\dfrac{\partial \mathcal{J}}{\partial (t,u)} \doteq \dfrac{\partial (y,x)}{\partial (t,u)} = \begin{pmatrix}
1 & 0 \\
L^1 & X^1 
\end{pmatrix} .
\end{equation}
From this we are able to compute the derivatives $\partial_y$ and $\partial_x$ in terms of $L$ and $\underline{\check{L}}$ as follows:
\begin{equation*}
\frac{\partial}{\partial y} = \frac{\partial}{\partial t} - \frac{L^1}{X^1} \frac{\partial}{\partial u} , \quad \quad \frac{\partial}{\partial x} = \frac{1}{X^1} \frac{\partial}{\partial u} ,
\end{equation*}
and since
$$ L = \frac{\partial}{\partial t} , \quad X = \frac{\partial}{\partial u}  , $$
we finally have that
\begin{equation}\label{xytotu}
\frac{\partial}{\partial y} = \left( 1 +\frac{\mu L^1}{2X^1} \right) L - \frac{L^1}{2X^1} \underline{\check{L}} , \quad \quad \frac{\partial}{\partial x} = - \frac{\mu}{2X^1} L + \frac{1}{2X^1} \underline{\check{L}} . 
\end{equation}

\subsection{The equations in dynamical coordinates}
We want to express the system of equations \eqref{eq:Phi} with respect to the $(L,\underline{\check{L}} )$ vector fields, where
$$ \underline{\check{L}} \doteq \mu \underline{L} , $$
and we note that
$$ \frac{\partial}{\partial y} = L , \quad \frac{\partial}{\partial x} = - \left( \frac{2L^1 + \mu}{2X^1} \right) L + \frac{1}{2X^1} \underline{\check{L}} . $$

For the nonlinearities we not that we have that:
$$ \mathcal{N}_1 ( \Phi_1 , \Phi_2 ) = H^{\nu \alpha \beta} [ ( \partial_{\beta} \Phi_{\alpha+1} ) ( \partial_{\nu} \Phi_1 ) - ( \partial_{\nu} \Phi_{\alpha + 1} ) ( \partial_{\beta} \Phi_1 ) ] - g^{\alpha \beta} \Omega^{\lambda} ( \partial_{\alpha} \Phi_{\lambda} ) ( \partial_{\beta} \Phi_1 ) , $$
$$ \mathcal{N}_2 ( \Phi_1 , \Phi_2 ) = H^{\nu \alpha \beta} [ ( \partial_{\beta} \Phi_{\alpha+1} ) ( \partial_{\nu} \Phi_1 ) - ( \partial_{\nu} \Phi_{\alpha + 1} ) ( \partial_{\beta} \Phi_2 ) ] - g^{\alpha \beta} \Omega^{\lambda} ( \partial_{\alpha} \Phi_{\lambda} ) ( \partial_{\beta} \Phi_2 ) , $$
where
$$ H^{\nu \alpha \beta} \doteq g^{\alpha \alpha_0 } g^{\beta \beta_0} H^{\nu}_{\alpha_0 \beta_0} , $$
$$ \Omega^{\lambda} = \frac{1}{\sqrt{-\det (g)}} \frac{\partial ( \sqrt{-\det(g)} )}{\partial \Phi_{\lambda}} , $$
for $\nu , \alpha , \alpha_0 , \beta , \beta_0 \in \{0,1\}$ and $\lambda \in \{1,2\}$.

Let
\begin{equation}\label{def:omega1}
 \omega_1 =- \frac{1}{2} H^{\gamma}_{\alpha \beta} L^{\alpha} L^{\beta} L_{\gamma} \check{X}^{\lambda} ( X \Phi_{\lambda}) , 
 \end{equation}
and
\begin{equation}\label{def:omega2}
 \omega_2 = - \frac{1}{2} H^{\gamma}_{\alpha \beta} L^{\alpha} L^{\beta} L_{\gamma}  \check{X}^{\lambda} ( L \Phi_{\lambda}) - \frac{1}{2} H^{\gamma}_{\alpha \beta} L^{\alpha} L^{\beta} \check{X}_{\gamma}  \check{X}^{\lambda} ( L \Phi_{\lambda}) . 
 \end{equation}
According to Proposition 4.3.1 of \cite{speck-book}, we have that
\begin{align*}
& \mu \Box_{g ( \Phi )} \Phi_{\kappa} = - L \check{\underline{L}} \Phi_{\kappa} ,
\\ & \mu \Box_{g (\Phi )} \Phi_{\kappa} = - \check{\underline{L}} L \Phi_{\kappa} - ( \omega_1 + \mu \omega_2 ) L \Phi_{\kappa}  ,
\end{align*}
for $\kappa \in \{1,2\}$. Hence, we will work with the following system of equations (derived from \eqref{eq:Phi}):
\begin{equation}\label{eq:Phi_null}
\begin{cases}
L \check{\underline{L}} \Phi_1  = \check{\mathcal{N}}_1 ( \Phi , L \Phi , \check{\underline{L}} \Phi ) , \\
\check{\underline{L}} L \Phi_1 + ( \omega_1 + \mu \omega_2 ) L \Phi_1  = \check{\mathcal{N}}_1 ( \Phi , L \Phi , \check{\underline{L}} \Phi ) , \\
 L \check{\underline{L}} \Phi_2   = \check{\mathcal{N}}_2 ( \Phi , L \Phi , \check{\underline{L}} \Phi ) , \\
\check{\underline{L}} L \Phi_2+ ( \omega_1 + \mu \omega_2 ) L \Phi_2 = \check{\mathcal{N}}_2 ( \Phi , L \Phi , \check{\underline{L}} \Phi ) ,
\end{cases}
\end{equation}
where
$$ \check{\mathcal{N}}_{\kappa} ( \Phi , L \Phi , \check{\underline{L}} \Phi ) = - \mu \mathcal{N}_{\kappa} ( \Phi_1 , \Phi_2 ) \mbox{  for $\kappa \in \{1,2\}$,} $$
and we note that both nonlinearities can be written only in terms of the $L$ and $\check{\underline{L}}$ of $\Phi_1$ and $\Phi_2$ by our previous computations.

Moreover, according to Lemma 3.2.1 from \cite{speck-book} we have the following equation for the inverse foliation density:
\begin{equation}\label{eq:mu}
L \mu = \omega_1 + \mu \omega_2 .
\end{equation}
The above equation involves not only the $H$ terms (that have been computed previously) and $L^0$ and $L^1$. By definition $L^0 = 1$, while for $L^1$ we have the following equation from Proposition 3.4.1 of \cite{speck-book}:
\begin{equation}\label{eq:l1}
L L^1 = \frac{1}{2} H_{\alpha \beta}^{\lambda} L^{\alpha} L^{\beta} ( L \Phi_{\lambda} ) \cdot \left( -L^1 - \frac{\Phi_1 \cdot \Phi_2}{c^2 - \Phi_1^2} \right) .
\end{equation}

\section{The main theorem}

In this section we state the main results that we are going to prove in this paper, and we also describe the method of proof. 

\begin{definition}
Let $T^*_U$ be the \textit{classical lifespan} of a solution of the system of equations  \eqref{eq:Phi} with data corresponding to the ones described in section \ref{data} which is the supremum of all $t > 0$ where $\Phi$ remains in $C^2 \times C^2$ in 
$$ \mathcal{R}_{T,U} \doteq \{ ( t,u ) \quad | \quad t \in [0,T) , u \in ( -\infty , U ] \} . $$
\end{definition}

Our main result gives an estimate on the size of $T*_U$ in terms of the sizes of $\Phi$ and $\check{\underline{L}} \Phi$ initially, and also gives an estimate on the behaviour of $\mu$ in $\mathcal{R}_{T,U}$.

\begin{theorem}\label{thm:main}
Let $\left( \left. ( \Phi_1 , \frac{\partial \Phi_1}{\partial n} ) \right|_{\Sigma} ,  \left. ( \Phi_2 , \frac{\partial \Phi_2}{\partial n} ) \right|_{\Sigma} \right)$ be an initial data set for the system of equations \eqref{eq:Phi_null} such that
\begin{equation}\label{data1}
\| \Phi_1 \|_{L^{\infty}} \leq \epsilon, \quad \| L \Phi_1 \|_{L^{\infty}} \leq \epsilon, \quad \| \underline{\check{L}} \Phi_1 \|_{L^{\infty}} \leq \delta , 
\end{equation}
\begin{equation}\label{data12}
\| L \Phi_2 \|_{L^{\infty}} \leq \epsilon, \quad \| \underline{\check{L}} \Phi_2 \|_{L^{\infty}} \leq \delta ,
\end{equation}
\begin{equation}\label{data3}
\sup_{\Sigma} | \Phi_2 | = K , \quad \inf_{\Sigma} | \Phi_2 | \geq K - \epsilon ,
\end{equation}
\begin{equation}\label{data4}
\| \Phi_1 \|_{L^2 ( \Sigma \cap \{ u \leq U_c \} )} \leq \epsilon , \quad \| \Phi_2 \|_{L^2 ( \Sigma \cap \{ u \leq U_c \} )} \leq \epsilon , 
\end{equation}
\begin{equation}\label{data5}
 \| L \Phi_1 \|_{L^2 ( \Sigma \cap \{ u \leq U_c \} )} \leq \epsilon , \quad \| L \Phi_2 \|_{L^2 ( \Sigma \cap \{ u \leq U_c \} )} \leq \epsilon ,
\end{equation}
for $\epsilon , \delta > 0$, $\epsilon$ small enough with respect to $\delta$, for some $U_c < \infty$, and $K $ and $\epsilon$ additionally satisfying  the following:
$$ \mbox{for  } \left. q \right|_{\Sigma} \doteq q_{in} \mbox{  we have that  } q_{in} > q_{cr} . $$

Then the classical lifespan time which is the maximal time (in $(t,u)$ coordinates) where the solution exists in $\cup_{t' \in [0,t)} \Sigma_{t'}$ and is $C^2$, is given by
$$ T^{*} = \sup\left\{t> 0 \quad | \quad \mu > 0 \mbox{  and  } \frac{q}{c} > 1 \right\} , $$
and as $\frac{q}{c} > 1$ is equivalent to $q > q_{cr}$, $T^*$ is also given by
$$ T^{*} = \sup\left\{t> 0 \quad | \quad \mu > 0 \mbox{  and  } q > q_{cr} \right\} . $$
For
$$ \delta_1 \doteq \sup_{\Sigma} [ \omega_1 ]_{-} \mbox{  and  } \delta_2 \doteq \sup_{\Sigma} [ \Phi_2 \cdot L \Phi_2 ]_{-} , $$
where $\omega_1$ was defined in \eqref{def:omega1} as
$$ \omega_1 =- \frac{1}{2} H^{\gamma}_{\alpha \beta} L^{\alpha} L^{\beta} L_{\gamma} X^{\lambda} ( \check{X} \Phi_{\lambda}) , $$
and where for a number $p$ we have that  $[p]_{-} \doteq |p|$  if $p < 0$ and  $[p]_{-} = 0$ otherwise, we have that
$$ 0 < T^* \leq \min \left( \frac{2}{\delta_1} , \frac{\min_u q^2(0,u) - q_{cr}^2}{2\delta_2} \right) , $$
and in particular we have that
\begin{equation}\label{vanishing}
\begin{split}
& T^* = \left(1+ O (\epsilon )\right)  \cdot  \frac{1}{\delta_1} \quad \mbox{    if $\frac{1}{\delta_1} <  \frac{\min_u q^2(0,u) - q_{cr}^2}{4\delta_2}$, or } 
\\ & T^* = \left(1+ O (\epsilon^2 )\right)  \cdot \frac{\min_u q^2(0,u) - q_{cr}^2}{4\delta_2}  \quad \mbox{    if $\frac{1}{\delta_1} >  \frac{\min_u q^2(0,u) - q_{cr}^2}{4\delta_2}$}.
\end{split} 
\end{equation}

For $(t,u) \in [0, T^{*} )$ we have the following estimates:
\begin{equation}\label{l1}
\| L \Phi_1 \|_{L^{\infty} (\Sigma_t )} + \| L \Phi_2 \|_{L^{\infty} (\Sigma_t )} \leq \| L \Phi_1 \|_{L^{\infty} (\Sigma )} + \| L \Phi_2 \|_{L^{\infty} (\Sigma )} + O ( \epsilon) ,
\end{equation}
\begin{equation}\label{l2}
\| \underline{\check{L}} \Phi_1 \|_{L^{\infty} (\Sigma_t )} + \| \underline{\check{L}} \Phi_2 \|_{L^{\infty} (\Sigma_t )} \leq \| \underline{\check{L}} \Phi_1 \|_{L^{\infty} (\Sigma )} + \| \underline{\check{L}} \Phi_2 \|_{L^{\infty} (\Sigma )} + O ( \epsilon) ,
\end{equation}
\begin{equation}\label{l3}
\mu (t,u) = c + C K ( \underline{\check{L}} \Phi_1 + \underline{\check{L}} \Phi_2 ) t + O ( \epsilon ) ,
\end{equation}
for some constants $c$, $C.$. Note that the above estimates imply that
$$ \| \Phi_1 \|_{L^{\infty} ( \Sigma_t)} \leq \bar{C} \epsilon , $$
$$  \| L \Phi_1 \|_{L^{\infty} ( \Sigma_t)} + \| L \Phi_2 \|_{L^{\infty} ( \Sigma_t)} \leq \bar{C} \epsilon , $$
$$ \| \underline{\check{L}} \Phi_1 \|_{L^{\infty} ( \Sigma_t)} + \| \underline{\check{L}} \Phi_2 \|_{L^{\infty} ( \Sigma_t)} \leq \bar{C} \delta , $$
and
$$ \| \Phi_2 \|_{L^{\infty} ( \Sigma_t)} \leq \bar{C} K , $$
for some constant $\bar{C}$.

In particular, $\underline{L} \Phi$ blows up along
$$ \{ ( T^{*} , u ) \quad | \quad \mu ( T^{*} , u ) = 0 \} \subseteq \Sigma_{T^{*}} \mbox{  if $\frac{1}{\delta_1} <  \frac{\min_u q^2(0,u) - q_{cr}^2}{4\delta_2}$,} $$
while a sonic line forms for 
$$\left\{ (T^{*} , u ) \quad | \quad \frac{q(T^{*} , u )}{c (T^{*} , u )} = 1 \right\} \subseteq \Sigma_{T^{*}}\mbox{  or  }  \left\{ (T^{*} , u ) \quad | \quad q(T^{*} , u ) = q_{cr} \right\} \subseteq \Sigma_{T^{*}} \mbox{ if $\frac{1}{\delta_1} >  \frac{\min_u q^2(0,u) - q_{cr}^2}{4\delta_2}$.} $$

Finally for all $t$ such that $0 \leq t < T^{*}$ we have that
$$ \lim_{u \rightarrow -\infty} \left. \Phi_1 (x) \right|_{\Sigma_t} = 0 , \quad   \lim_{u \rightarrow -\infty} \left. \Phi_2 (x) \right|_{\Sigma_t} = K . $$
\end{theorem}

It is worth noticing that the boundary of the spacetime region of existence of the solution can consist not only of points where a shock forms, i.e. where $\mu = 0$, but also of points where a sonic line forms, i.e. where $\frac{q}{c} = 1$. Moreover, the last part of the Theorem essentially states that the Mach number at infinity is fixed throughout the region of existence.

Note that for data
$$ \left. \Phi_1 \right|_{\Sigma} = 0, \quad \left. \frac{\partial \Phi_1}{\partial n} \right|_{\Sigma} \sim \epsilon , \quad \left. \Phi_2 \right|_{\Sigma} \sim q_{cr} , \quad L \Phi_2 \sim \epsilon , \quad  \underline{\check{L}} \Phi_2 \sim \delta , $$
we expect that in the region where the solution exists we have that
$$ c^2 \sim q_{cr}^2 , \quad c^2 - \Phi_1^2 \sim q_{cr}^2 , \quad c^2 - \Phi_1^2 - \Phi_2^2 \sim -\epsilon^2 .$$ 

Note then that for the nonlinearities $\mathcal{N}_1$ and $\mathcal{N}_2$ -- taking into consideration the size of $K$ and $\epsilon$ -- we have schematically that
$$ \mathcal{N}_1 \sim \mathcal{N}_2 \sim \sum_{i,j,k \in \{1,2\}} ( 1+ \Phi_1 \Phi_2 ) \left[ \Phi_i (L \Phi_j ) (L \Phi_k ) + \Phi_i ( L \Phi_j ) (\underline{L} \Phi_k )  \right]  , $$
and as $\epsilon$ is taken to be small, we essentially have -- once again schematically -- that
$$  \mathcal{N}_1 \sim \mathcal{N}_2 \sim \sum_{i,j,k \in \{1,2\}}  \Phi_i (L \Phi_j ) (L \Phi_k )  + \Phi_i ( L \Phi_j ) ( \underline{L} \Phi_k )  . $$

Writing the system of equation \eqref{eq:Phi} using the vector fields $(L , \underline{\check{L}} )$ schematically we obtain the following system
\begin{equation}\label{eq:schem}
\begin{cases}
L \underline{\check{L}} \Phi_1 \sim \sum_{i,j,k \in \{1,2\}} \left[ \mu \Phi_i ( L \Phi_j ) ( L \Phi_k ) + \Phi_i ( L \Phi_j ) ( \underline{\check{L}} \Phi_k ) \right]  , \\
\underline{\check{L}} L \Phi_1 \sim \sum_{i,j,k \in \{1,2\}} \left[ \mu \Phi_i ( L \Phi_j ) ( L \Phi_k ) + \Phi_i ( L \Phi_j ) ( \underline{\check{L}} \Phi_k ) \right] , \\ 
L \underline{\check{L}} \Phi_2 \sim \sum_{i,j,k \in \{1,2\}} \left[ \mu \Phi_i ( L \Phi_j ) ( L \Phi_k ) + \Phi_i ( L \Phi_j ) ( \underline{\check{L}} \Phi_k ) \right]  , \\
\underline{\check{L}} L \Phi_2 \sim \sum_{i,j,k \in \{1,2\}} \left[ \mu \Phi_i ( L \Phi_j ) ( L \Phi_k ) + \Phi_i ( L \Phi_j ) ( \underline{\check{L}} \Phi_k ) \right] .
\end{cases}
\end{equation}

Note that the inverse foliation density $\mu$ satisfies a schematic equation of the following form (see Lemma A.2.1 of \cite{speck-book}):
\begin{equation}\label{eq:mu_schem}
L \mu \sim \sum_{i,j \in \{1,2\}} \Phi_i ( \mu L \Phi_j + \underline{\check{L}} \Phi_j ) .
\end{equation}

We will argue through a bootstrap argument. We will assume that:
$$ \| \Phi_1 \|_{L^{\infty} (\Sigma_t )} \leq c_1 \epsilon , \quad \| L \Phi_1 \|_{L^{\infty} (\Sigma_t )} \leq c_2 \epsilon , \quad \| L \Phi_2 \|_{L^{\infty} (\Sigma_t )} \leq c_3 \epsilon  , $$
for some $c_1, c_2, c_3 > 0$, 
$$ \| \Phi_2 \|_{L^{\infty} (\Sigma_t )} \leq d_1 K , \quad \| \underline{\check{L}} \Phi_1 \|_{L^{\infty} (\Sigma_t )} \leq d_2 \delta  , \quad \| \underline{\check{L}} \Phi_2 \|_{L^{\infty} (\Sigma_t )} \leq d_3 \delta   , $$
for some $d_1 , d_2 , d_3 > 0$, and
\begin{equation*}
\begin{split}
&  \| \Phi_1 \|_{L^2 ( \Sigma \cap \{ u \leq U_c \} )} \leq e_1 \epsilon , \quad \| \Phi_2 - K \|_{L^2 ( \Sigma \cap \{ u \leq U_c \} )} \leq e_2 \epsilon ,  \\ & \| L \Phi_1 \|_{L^2 ( \Sigma \cap \{ u \leq U_c \} )} \leq e_3 \epsilon , \quad \| L \Phi_2 \|_{L^2 ( \Sigma \cap \{ u \leq U_c \} )} \leq e_4 \epsilon , 
\end{split}
\end{equation*} 
for some $e_1 , e_2 , e_3 , e_4 > 0$. The goal is to improve the $c_i$'s, the $d_i$'s, and the $e_i$'s.

We start by constructing the solution close to spacelike infinity. Note that an outcome of the estimates we will establish is that the Mach number at infinity remains fixed (which can be seen as a boundary condition). We will use a basic energy estimate. We will show the following estimate:
\begin{align*}
\int_{\Sigma_t} & \Big[ (1+\mu ) \mu  ( L \Phi_1 )^2+  \mu ( L \Phi_1 ) ( \check{X} \Phi_1 ) +  ( \check{X} \Phi_1 )^2 \Big] \, du \\ & + \int_{\Sigma_t} \Big[ (1+\mu ) \mu  ( L (\Phi_2 - K ) )^2+  \mu ( L (\Phi_2 - K ) ) ( \check{X}( \Phi_2 - K ) ) +  ( \check{X}(  \Phi_2 - K ) )^2 \Big] \, du\\ \leq & \int_{\Sigma_0} \Big[ (1+\mu ) \mu  ( L \Phi_1 )^2+  \mu ( L \Phi_1 ) ( \check{X} \Phi_1 ) +  ( \check{X} \Phi_1 )^2 \Big] \, du \\ &  + \int_{\Sigma_0} \Big[ (1+\mu ) \mu  ( L ( \Phi_2 - K ) )^2+  \mu ( L ( \Phi_2 - K ) ) ( \check{X} ( \Phi_2 - K ) ) +  ( \check{X} (\Phi_2 - K ) )^2 \Big] \, du \\ & + O (\epsilon^3 ) ,
\end{align*}
using the last set of the aforementioned bootstrap assumptions, which will allows us to conclude that
\begin{align*}
\int_{\Sigma_t} & \Big[ (1+\mu ) \mu  ( L \Phi_1 )^2+  \mu ( L \Phi_1 ) ( \check{X} \Phi_1 ) +  ( \check{X} \Phi_1 )^2 \Big] \, du \\ & + \int_{\Sigma_t} \Big[ (1+\mu ) \mu  ( L (\Phi_2 - K ) )^2+  \mu ( L (\Phi_2 - K ) ) ( \check{X}( \Phi_2 - K ) ) +  ( \check{X}(  \Phi_2 - K ) )^2 \Big] \, du < \infty , 
\end{align*}
for all $t \in [ 0 , T^{*} )$. The result will follow by the coercive estimate
\begin{align*}
 \int_{\Sigma_t}\Big[ &  \Phi_1^2 + ( \Phi_2 - K )^2 \Big] \, du \leq  \int_{\Sigma_t}  \Big[ (1+\mu ) \mu  ( L \Phi_1 )^2+  \mu ( L \Phi_1 ) ( \check{X} \Phi_1 ) +  ( \check{X} \Phi_1 )^2 \Big] \, du \\ & + \int_{\Sigma_t} \Big[ (1+\mu ) \mu  ( L (\Phi_2 - K ) )^2+  \mu ( L (\Phi_2 - K ) ) ( \check{X}( \Phi_2 - K ) ) +  ( \check{X}(  \Phi_2 - K ) )^2 \Big] \, du ,
\end{align*} 
again for $t \in [0, T^{*} )$. We note that we need to work with a slightly different system compared to \eqref{eq:Phi}, one satisfied by $\Phi_1$ and $\Phi_2 - K$.

Next we will construct the solution in the interior. We will use equations \eqref{eq:schem} and show first the estimates \eqref{l1}, \eqref{l2}, and \eqref{l3}. We start with $L \Phi$, and use the second and fourth equations of \eqref{eq:schem} and we have that for $k \in \{1,2\}$:
$$ | \underline{\check{L}} L \Phi_k | (t,u) \leq \sum_{i=1}^2 K \delta | L \Phi_i | (t,u) , $$
and considering $u$ in a bounded (potentially small) interval, we can apply a Gr\"{o}nwall estimate using the previous inequality (integrating along the integral curves of $\underline{\check{L}}$), and we have that for $k \in \{1,2\}$:
$$ | L \Phi_k | (t , u) \leq | L \Phi_k | (t, u_0) + C_1 K \delta \epsilon + O (\epsilon^2 ) . $$ 
Integrating then along the integral curves of $L$ we have that
$$ | \Phi_1 | (t) \leq | \Phi_1 | (0) + C_2 K \delta \epsilon + O (\epsilon^2 ) \leq K + C_2 K \delta \epsilon + O ( \epsilon^2 ), $$
$$ | \Phi_2 | (t) \leq | \Phi_2 | (0) + C_3 K \delta \epsilon + O (\epsilon^2 ) \leq \epsilon + C_3 K \delta \epsilon + O ( \epsilon^2 ) . $$

Having established the previous estimates we can now turn to $\underline{\check{L}}$ using the first and third equations of \eqref{eq:schem}, and we have that for $k \in \{1,2\}$:
$$ | \underline{\check{L}} \Phi_k | (t) \leq | \underline{\check{L}} \Phi_k | (0) + C_4 K \delta \epsilon + O ( \epsilon^2 ) \leq \delta + C_4 K \delta \epsilon + O ( \epsilon^2 ) . $$

Concerning shock formation, for $\mu$ we use equation \eqref{eq:mu} and we have that
$$ L \mu \sim K ( \underline{\check{L}} \Phi_1 + \underline{\check{L}} \Phi_2 ) + O ( \epsilon ) . $$
After integrating along the integral curves of $L$ we have that
$$ \mu (t) \sim \mu (0) + K  ( \underline{\check{L}} \Phi_1 + \underline{\check{L}} \Phi_2 ) t + O ( \epsilon ) , $$
and the definition of $\delta_0$ as the quantity controlling the blow up time follows from the previous calculation. Note that in the area close to spacelike infinity it can be shown through the same method that $\mu$ remains positive and bounded away from zero, so a shock can form only in the interior region.

Finally, for the formation of the sonic line we consider the quantity $L ( \Phi_1^2 + \Phi_2^2 )$  that we integrate in $t$. We use again the established estimates for $\Phi_i$, $L \Phi_i$, $i \in \{1,2\}$ in order to compute the estimated time where $\Phi_1^2 + \Phi_2^2 = q_{cr}^2$. Note again that the previously established estimates for $\Phi_i$, $L \Phi_i$, $i \in \{1,2\}$ guarantee that close to spacelike infinity we have that $\Phi_1^2 + \Phi_2^2$ remains bounded away from $q_{cr}^2$, so a sonic line can form only in the interior region.

\begin{remark}
We note that $\epsilon$ is chosen to be small enough not only in order to be able to close the bootstrap assumptions, but also in order to ensure that equation \eqref{eq:phi} remains hyperbolic throughout the region where the solution exists, i.e. that
$$ q > q_{cr} \mbox{  for all $0 \leq t < T^{*}$.} $$
\end{remark}

\section{Proof of the main theorem}
In this section we prove the main Theorem \ref{thm:main}. We start by constructing the set of initial data that we consider. The proof follows by a standard bootstrap argument that occupies the largest part of this section. We will state certain bootstrap assumptions that will imply certain estimates on the nonlinearities. Using these we will be able to close the bootstrap argument on an energy and pointwise level.

 As we need to work with energy spaces we define the following renormalized quantities:
$$ \Psi_1 \doteq \Phi_1 , \quad \quad \Psi_2 \doteq \Phi_2 - K.$$
The new quantity $\Psi = ( \Psi_1 , \Psi_2 )$ satisfies the following system of equations
\begin{equation}\label{eq:Psi}
\begin{cases}
\Box_{\mathbf{g}( \Psi )} \Psi_1 = \mathbf{N}_1 ( \Psi_1 , \Psi_2 ) , \\
\Box_{\mathbf{g} (\Psi )} \Psi_2 = \mathbf{N}_2 (\Psi_1 , \Psi_2 ) ,
\end{cases}
\end{equation}
where
\begin{equation*}
\mathbf{g}^{-1} ( \Psi ) =
\begin{bmatrix}
-1 & -\frac{ \Psi_1 \cdot ( \Psi_2 + K )}{\mathbf{c}^2 - \Psi_1^2} \\
-\frac{ \Psi_1 \cdot ( \Psi_2 +K )}{\mathbf{c}^2 - \Psi_1^2} & - \frac{ \mathbf{c}^2- ( \Psi_2 + K )^2}{\mathbf{c}^2 - \Psi_1^2} 
\end{bmatrix} ,
\end{equation*}
where
$$ \mathbf{c} ( \Psi ) = c ( \Psi_1^2 + ( \Psi_2 + q_{cr} )^2 ) , $$
and where
\begin{align*}
\mathbf{N}_1 ( \Psi_1 , \Psi_2) & = \frac{\partial}{\partial x} \left( \frac{\Psi_1 \cdot ( \Psi_2 + K )}{\mathbf{c}^2 - \Psi_1^2} \right) \partial_y \Psi_1 - \frac{\partial}{\partial y} \left( \frac{\Psi_1 \cdot ( \Psi_2 + K )}{\mathbf{c}^2 - \Psi_1^2} \right) \partial_x \Psi_1 \\ & + \frac{\partial}{\partial y} \left( \frac{\mathbf{c}^2 - ( \Psi_2 + K )^2}{\mathbf{c}^2 - \Psi_1^2} \right) \partial_x \Psi_2 - \frac{\partial}{\partial x} \left( \frac{\mathbf{c}^2 - ( \Psi_2 + K)^2}{\mathbf{c}^2 - \Psi_1^2} \right) \partial_x \Psi_1 \\ & + \frac{1}{\sqrt{-\det (\mathbf{g})}} [ \partial_{\alpha}  ( \sqrt{-\det (\mathbf{g})} ) ] \cdot  \mathbf{g}^{\alpha \beta} \partial_{\beta} \Psi_1 ,
\end{align*}
\begin{align*}
\mathbf{N}_2 ( \Psi_1 , \Psi_2 ) & = \frac{\partial}{\partial y} \left( \frac{\Psi_1 \cdot (\Psi_2 + K )}{ \mathbf{c}^2 - \Psi_1^2} \right) \partial_x \Psi_2 - \frac{\partial}{\partial x} \left( \frac{\Psi_1 \cdot ( \Psi_2 + K )}{ \mathbf{c}^2 - \Psi_1^2} \right) \partial_y \Psi_2 \\ & + \frac{1}{\sqrt{-\det (\mathbf{g})}} [ \partial_{\alpha}  ( \sqrt{-\det (\mathbf{g})} ) ] \cdot  \mathbf{g}^{\alpha \beta} \partial_{\beta} \Psi_2 .
\end{align*}

\subsection{Construction of the initial data set}
We can choose any smooth functions $(f_1 , f_2 )$ in one variable for 
$$ f_1 = \left. \Phi_1 \right|_{\Sigma} , \quad \quad f_2 = \left. \Phi_2 \right|_{\Sigma} , $$
with the appropriate size restrictions in $L^{\infty}$ and $L^2$. If the size restrictions of the derivatives of $(f_1 , f_2 )$ are not satisfied, then we can argue by rescaling in order to get the desired bounds. 

On the other hand, a consequence of \eqref{def:lmu} is that
$$ \left. \mu \right|_{\Sigma}= \frac{1}{\sqrt{(\underline{g}^{-1})^{00}}} , $$
for $\underline{g} = \Pi g$, $\Pi$ the projection to $\Sigma$ (for a proof see Lemma 7.2 of \cite{jaredgustavjonwillie}). Initially at $\Sigma$ this gives us that $\inf_{\Sigma} \mu = 1$.

\subsection{The bootstrap assumptions}
Let us make the following bootstrap assumptions: for some constant $\mathcal{C}$ we have that
\begin{equation}\label{boot1}
\| \Phi_1 \|_{L^{\infty} (\Sigma_t )} \leq \mathcal{C} \epsilon ,
\end{equation}
\begin{equation}\label{boot2}
\| L \Phi_1 \|_{L^{\infty} (\Sigma_t )} + \|L  \Phi_2 \|_{L^{\infty} (\Sigma_t )} \leq \mathcal{C} \epsilon ,
\end{equation}
\begin{equation}\label{boot3}
\| \check{\underline{L}} \Phi_1 \|_{L^{\infty} (\Sigma_t )}+ \| \check{\underline{L}} \Phi_2 \|_{L^{\infty} (\Sigma_t )} \leq \mathcal{C} \delta ,
\end{equation}
\begin{equation}\label{boot4}
\| \Phi_2 \|_{L^{\infty} (\Sigma_t )} \leq \mathcal{C} K ,
\end{equation}
for $\epsilon$, $\delta$ and $K$ as in the assumptions and for $t \in [ 0 , T^0 )$ for $T^0$ to be a bootstrap time such that
$$ 0 < T^0 < \min \left( \frac{2}{\delta_1} , \frac{\min_u q^2(0,u) - q_{cr}^2}{2\delta_2} \right) . $$

Moreover we also assume that again for $t \in [ 0 , T^0 )$ we have that
\begin{equation}\label{boot5}
\| \Phi_1 \|_{L^{\infty} (\Sigma_t \cap \{ u \leq U_c \} )}+\| \Phi_2 - K \|_{L^{\infty} (\Sigma_t \cap \{ u \leq U_c \} )} \leq \underline{\mathcal{C}} \epsilon ,
\end{equation}
\begin{equation}\label{boot6}
\| L \Phi_1 \|_{L^{\infty} (\Sigma_t \cap \{ u \leq U_c \})} + \|L  \Phi_2 \|_{L^{\infty} (\Sigma_t \cap \{ u \leq U_c \})} \leq \underline{\mathcal{C}} \epsilon ,
\end{equation}
\begin{equation}\label{boot7}
\| \check{\underline{L}} \Phi_1 \|_{L^{\infty} (\Sigma_t \cap \{ u \leq U_c \})}+ \| \check{\underline{L}} \Phi_2 \|_{L^{\infty} (\Sigma_t \cap \{ u \leq U_c \} )} \leq \underline{\mathcal{C}} \epsilon .
\end{equation}

Under the bootstrap assumptions \eqref{boot1}, \eqref{boot2}, \eqref{boot3} and \eqref{boot4} we will show that:
\begin{equation}\label{boot:non1}
 | \check{\mathcal{N}}_1 ( \Phi , L \Phi , \check{\underline{L}} \Phi ) | + | \check{\mathcal{N}}_2 ( \Phi , L \Phi , \check{\underline{L}} \Phi ) | + ( | \omega_1 | + \mu | \omega_2 | ) | L \Phi_{\kappa} | \leq \sum_{i=1}^2 \mathcal{C}' K \delta  | L \Phi_i | , 
 \end{equation}
and 
for $\kappa \in \{1,2\}$ and for all $(t,u) \in [0, T^0 ) \times (-\infty , U]$. Note that for $(t,u) \in [0, T^0 ) \times (-\infty , U_0]$ we do not need the bootstrap assumptions but we can just use the estimates of the previous section. This also implies that:
\begin{equation}\label{boot:non2}
 | \omega_1 | + \mu | \omega_2 | \leq \mathcal{C}'' K \delta ,
\end{equation}
\begin{equation}\label{boot:non3}
 | \check{\mathcal{N}}_1 ( \Phi , L \Phi , \check{\underline{L}} \Phi ) | + | \check{\mathcal{N}}_2 ( \Phi , L \Phi , \check{\underline{L}} \Phi ) | + ( | \omega_1 | + \mu | \omega_2 | ) | L \Phi_{\kappa} | \leq  \mathcal{C}''' K \delta \epsilon, 
 \end{equation}
again for $\kappa \in \{1,2\}$ and for all $(t,u) \in [0, T^0 ) \times (-\infty , U]$, and it also implies that:
\begin{equation}\label{boot:non4}
 | \mathbf{N}_1 ( \Psi_1 , \Psi_2 ) | + | \mathbf{N}_2 ( \Psi_1 , \Psi_2 ) |  \leq   \sum_{i=1}^2\mathcal{C}'''' \epsilon^2 | L \Psi_i | ,
 \end{equation}
for all $(t,u) \in [0, T^0 ) \times (-\infty , U_c ]$.

We start by showing \eqref{boot:non1}, examining first the terms $\mathcal{N}_1$ and $\mathcal{N}_2$. From Lemma A.3.1 of \cite{speck-book} we have that
$$  \mathcal{N}_1   = f_{11} \cdot ( L \Phi_ 1 ) ( L \Phi_2 ) + f_{12} (L \Phi_1 ) ( X \Phi_2 ) + f_{13} ( X \Phi_1 ) ( L \Phi_2 ) , $$
$$  \mathcal{N}_2   = f_{21} \cdot ( L \Phi_ 1 ) ( L \Phi_2 ) + f_{22} (L \Phi_1 ) ( X \Phi_2 ) + f_{23} ( X \Phi_1 ) ( L \Phi_2 ) , $$
where the $f_i$'s, $i \in \{1,2,3\}$, are smooth functions of the $\Phi_j$'s, $j \in \{1,2\}$ and the components $L^j$'s and $X^j$'s, $j \in \{1,2\}$. From the above estimate \eqref{boot:non1} for $\mathcal{N}_1$, $\mathcal{N}_2$, follows from the bootstrap assumptions \eqref{boot1}, \eqref{boot2}, \eqref{boot3} and \eqref{boot4} and the form of the nonlinearities, after showing again through a bootstrap argument that
\begin{equation}\label{boot:lx}
| L^i | \leq c_L, \quad | X^i | \leq c_X \mbox{  for $i \in \{1,2\}$}, 
\end{equation}
for some constants $c_L > 0$, $c_X > 0$. This can be done similarly to Lemma 11.9.1 from \cite{speck-book}, with obvious modifications for $u \in (-\infty , U_c]$. 

From the form of the $\mathcal{N}_i$'s, $i \in \{1,2\}$, we have that
$$ |f_{ij} | \leq \sum_{k=1}^2 c_f | \Phi_k |  \mbox{  for $i \in \{1,2,3\}$, $j \in \{1,2\}$, and a constant $c_f > 0$.} $$  
This completes the proof of estimate \eqref{boot:non1} for $\mathcal{N}_1$ and $\mathcal{N}_2$. Estimate \eqref{boot:non1} for the $\omega$'s follows from the form of the $H^{\nu}_{\alpha \beta}$'s, while estimate \eqref{boot:non3} and \eqref{boot:non4} are consequences of estimate \eqref{boot:non1}.

\subsection{Construction of the solution near spacelike infinity}

First we construct the solution close to spacelike infinity, that is for $|x| \rightarrow \infty$ along $\Sigma_t$ for any $t \in [0, T^0 )$ (the fact that this remains fixed will follow by the construction as the $\Sigma_t$ hypersurfaces depend on the solution $\Phi$), also showing that throughout the region of existence and uniqueness of the solution $\Phi$, the Mach speed at infinity will also remain fixed (so there is no need imposing a boundary condition). Eventually after closing the bootstrap argument we will have the same result for $t \in [0, T^* )$.

The estimates that we will establish will imply that
$$ \lim_{u \rightarrow -\infty} \sqrt{ \Psi_1^2 (t,u) + \Psi_2^2 (t,u) } = 0 , $$
ifor all $t \in [0, T^{*} )$, hence ensuring that the Mach speed at infinity remains fixed.

Recall that the energy-momentum tensor for the wave equation is given by:
$$ \mathbb{T}_{\alpha \beta} [ \psi ] \doteq \partial_{\alpha} \psi \cdot \partial_{\beta} \psi - \frac{1}{2} \mathbf{g}_{\alpha \beta} \mathbf{g}^{\sigma \tau} \partial_{\sigma} \psi \cdot \partial_{\tau} \psi \mbox{  for $\alpha , \beta , \sigma , \tau \in \{0,1\}$} . $$ 
For the system \eqref{eq:Psi} we have the following well-known computations:
$$ \partial_{\alpha} \mathbb{T}^{\alpha \beta} [ \Psi_1 ] = \mathbf{N}_1 ( \Psi_1 , \Psi_2 ) \cdot \partial^{\beta} \Psi_1 \mbox{  and  }  \partial_{\alpha} \mathbb{T}^{\alpha \beta} [ \Psi_2 ] = \mathbf{N}_2 ( \Psi_1 , \Psi_2 ) \cdot \partial^{\beta} \Psi_2 . $$
For a vector field $V$ we define its energy current as:
$$ \mathbb{J}^V_{\alpha} [ \psi ] \doteq \mathbb{T}_{\alpha \beta} [ \psi ] \cdot V^{\beta} . $$
For the system \eqref{eq:Psi} we have another well-known computation:
$$ \mathrm{Div} ( J^V [ \Psi_1 ] ) = \mathbf{N}_1 ( \Psi_1 , \Psi_2 ) \cdot V\Psi_1 + \mathbb{T}_{\alpha \beta} [ \Psi_1 ] \cdot \pi_V^{\alpha \beta} \mbox{  and  } $$ 
$$\mathrm{Div} ( J^V [ \Psi_2 ] ) = \mathbf{N}_2 ( \Psi_1 , \Psi_2 ) \cdot V\Psi_1 + \mathbb{T}_{\alpha \beta} [ \Psi_2 ] \cdot \pi_V^{\alpha \beta} , $$
where $\pi_V$ is the deformation tensor with respect to $V$ defined as:
$$ ( \pi_V )_{\alpha \beta} \doteq \partial_{\alpha} V_{\beta} + \partial_{\beta} V_{\alpha} . $$ 

Let us choose some $-\infty < U_0 < U_c < U$. An application of the divergence theorem gives the following identity (for a detailed proof see Lemma 3.1 of \cite{speck-notes}):
\begin{equation}\label{eq:div}
\begin{split}
\int_{\Sigma_T^{U_0}} \mu J^V_{\alpha} [ \psi ] \cdot n_{\Sigma_T}^{\alpha} \, du +&  \int_{\mathcal{N}_{U_0}^T}J^V_{\alpha} [ \psi ] \cdot L^{\alpha} \, du  =  \int_{\Sigma_0^{U_0}} \mu J^V_{\alpha} [ \psi ] \cdot n_{\Sigma} \, du \\  &-   \int_{\mathcal{R}_{T,U_0}} \mu \mathbb{T}_{\alpha \beta} [ \psi ] \cdot \pi_V^{\alpha \beta} \, dt du - \int_{\mathcal{R}_{T,U_0}} \mu \Box_{\mathbf{g}} \psi \cdot V \psi \, dt du ,
\end{split}
\end{equation}
for
$$ \Sigma_{T}^{U_0} \doteq \{ t = T , -\infty < u \leq U_0 \} , \quad \mathcal{N}_{U_0}^T = \{ u = U_0 , 0 \leq t < T \} , $$
for $\mathcal{R}_{T,U_0}$ as defined before, and for $n_{\Sigma}$ the normal to the surface $\Sigma$. We choose $T \leq T^0$.

Choosing 
$$ V = (1+\mu ) L + \check{\underline{L}} = (1+2\mu ) L + 2 \check{X} , $$
and using the aforementioned \eqref{eq:div} we have that
\begin{equation}\label{eq:enest}
\begin{split}
\int_{\Sigma_T^{U_0}} \Big[ \frac{1}{2} ( 1+2\mu ) & \mu ( L \psi )^2+ 2 \mu ( L \psi ) ( \check{X} \psi ) + 2 ( \check{X} \psi )^2 \Big] \, du +  \int_{\mathcal{N}_{U_0}^T} ( 1+ \mu ) ( L \psi )^2 \, du  \\ = &  \int_{\Sigma_0^{U_0}} \Big[ \frac{1}{2} ( 1+2\mu ) \mu ( L \psi )^2 + 2 \mu ( L \psi ) ( \check{X} \psi ) + 2 ( \check{X} \psi )^2 \Big] \, du \\  &-   \int_{\mathcal{R}_{T,U_0}} \mu \mathbb{T}_{\alpha \beta} [ \psi ] \cdot \pi_V^{\alpha \beta} \, dt du - \int_{\mathcal{R}_{T,U_0}} \mu \Box_{\mathbf{g}} \psi \cdot [ (1+2\mu ) (L \psi ) +2 (\check{X} \psi ) ] \, dt du .
\end{split}
\end{equation}
For a proof see Lemma 3.2 and Proposition 3.4 of \cite{speck-notes}. Applying the identity \eqref{eq:enest} to the solutions $\Psi_1$ and $\Psi_2$ of the system \eqref{eq:Psi} we have that:
\begin{equation}\label{eq:enestPsi}
\begin{split}
\int_{\Sigma_T^{U_0}} \Big[ \frac{1}{2} & ( 1+2\mu )  \mu ( L \Psi_1 )^2+ 2 \mu ( L \Psi_1 ) ( \check{X} \Psi_1 ) + 2 ( \check{X} \Psi_1 )^2 \Big] \, du +  \int_{\mathcal{N}_{U_0}^T} ( 1+ \mu ) ( L \Psi_1 )^2 \, du  \\ & +\int_{\Sigma_T^{U_0}} \Big[ \frac{1}{2} ( 1+2\mu )  \mu ( L \Psi_2 )^2+ 2 \mu ( L \Psi_2 ) ( \check{X} \Psi_2 ) + 2 ( \check{X} \Psi_2 )^2 \Big] \, du +  \int_{\mathcal{N}_{U_0}^T} ( 1+ \mu ) ( L \Psi_2 )^2 \, du  \\ = &  \int_{\Sigma_0^{U_0}} \Big[ \frac{1}{2} ( 1+2\mu ) \mu ( L \Psi_1 )^2 + 2 \mu ( L \Psi_1 ) ( \check{X} \Psi_1 ) + 2 ( \check{X} \Psi_1 )^2 \Big] \, du \\ & + \int_{\Sigma_0^{U_0}} \Big[ \frac{1}{2} ( 1+2\mu ) \mu ( L \Psi_2 )^2 + 2 \mu ( L \Psi_2 ) ( \check{X} \Psi_2 ) + 2 ( \check{X} \Psi_2 )^2 \Big] \, du\\  & -   \int_{\mathcal{R}_{T,U_0}} \mu \mathbb{T}_{\alpha \beta} [ \Psi_1 ] \cdot \pi_V^{\alpha \beta} \, dt du - \int_{\mathcal{R}_{T,U_0}} \mu \mathbb{N}_1 ( \Psi_1 , \Psi_2 ) \cdot [ (1+2\mu ) (L \Psi_1 ) +2 (\check{X} \Psi_1 ) ] \, dt du \\ & -   \int_{\mathcal{R}_{T,U_0}} \mu \mathbb{T}_{\alpha \beta} [ \Psi_2 ] \cdot \pi_V^{\alpha \beta} \, dt du - \int_{\mathcal{R}_{T,U_0}} \mu \mathbb{N}_2 ( \Psi_1 , \Psi_2 ) \cdot [ (1+2\mu ) (L \Psi_2 ) +2 (\check{X} \Psi_2 ) ] \, dt du .
\end{split}
\end{equation}
We note that according to Proposition 3.4 of \cite{speck-notes} we have that
\begin{equation}\label{eq:deformation}
 \mu \mathbb{T}_{\alpha \beta} [ \psi ] \cdot \pi_V^{\alpha \beta} = ( L \psi )^2 \Big[ -\frac{1}{2} ( L \mu ) + \mu ( \check{X} \mu ) \Big] . 
 \end{equation}

By Proposition 3.3 of \cite{speck-notes} and an application of the Cauchy-Schwarz inequality in the estimate \eqref{eq:enestPsi} we have that
\begin{equation}\label{eq:enestPsi1}
\begin{split}
\int_{\Sigma_T^{U_0}} \Big[ \frac{1}{2} & ( 1+2\mu )  \mu ( L \Psi_1 )^2 + 2 ( \check{X} \Psi_1 )^2 \Big] \, du +  \int_{\mathcal{N}_{U_0}^T} ( 1+ \mu ) ( L \Psi_1 )^2 \, du  \\ & +\int_{\Sigma_T^{U_0}} \Big[ \frac{1}{2} ( 1+2\mu )  \mu ( L \Psi_2 )^2 + 2 ( \check{X} \Psi_2 )^2 \Big] \, du +  \int_{\mathcal{N}_{U_0}^T} ( 1+ \mu ) ( L \Psi_2 )^2 \, du  \\ \leq &  \int_{\Sigma_0^{U_0}} \Big[ 2( 1+2\mu ) \mu ( L \Psi_1 )^2+ 4 ( \check{X} \Psi_1 )^2 \Big] \, du \\ & + \int_{\Sigma_0^{U_0}} \Big[2( 1+2\mu ) \mu ( L \Psi_2 )^2  + 4 ( \check{X} \Psi_2 )^2 \Big] \, du\\  & -   \int_{\mathcal{R}_{T,U_0}} \mu \mathbb{T}_{\alpha \beta} [ \Psi_1 ] \cdot \pi_V^{\alpha \beta} \, dt du - \int_{\mathcal{R}_{T,U_0}} \mu \mathbb{N}_1 ( \Psi_1 , \Psi_2 ) \cdot [ (1+2\mu ) (L \Psi_1 ) +2 (\check{X} \Psi_1 ) ] \, dt du \\ & -   \int_{\mathcal{R}_{T,U_0}} \mu \mathbb{T}_{\alpha \beta} [ \Psi_2 ] \cdot \pi_V^{\alpha \beta} \, dt du - \int_{\mathcal{R}_{T,U_0}} \mu \mathbb{N}_2 ( \Psi_1 , \Psi_2 ) \cdot [ (1+2\mu ) (L \Psi_2 ) +2 (\check{X} \Psi_2 ) ] \, dt du \\ \doteq & I + II + III + IV .
\end{split}
\end{equation}
Note that terms $I$ and $III$ are similar, while the same is true for terms $II$ and $IV$. We start with the last two, and by bootstrap assumption \eqref{boot:non4} and Cauchy-Schwarz we have that
$$ |II| + |IV| \leq 2 \mathcal{C}'''' \epsilon^2  \int_{\mathcal{R}_{T,U_0}}  (1+2\mu ) \mu [ (L \Psi_1 )^2 + (L \Psi_2 )^2 ]  +2 [ ( \check{X} \Psi_1 )^2 + (\check{X} \Psi_2 )^2 ] \, dt du . $$
Using the fact that $T \leq \min \left( \frac{2}{\delta_1} , \frac{\min_u q^2(0,u) - q_{cr}^2}{2\delta_2} \right)$ and that neither $\delta_1$ nor $\delta_2$ are taken to be small (while $\epsilon$ is), we have that
$$ |II| + |IV| \leq \bar{\mathcal{C}}'''' \epsilon^2 \sup_{ t \in [ 0 ,T ]} \int_{\Sigma_t \cap \{ u \leq U_0\}}    (1+2\mu ) \mu [ (L \Psi_1 )^2 + (L \Psi_2 )^2 ]  +2 [ ( \check{X} \Psi_1 )^2 + (\check{X} \Psi_2 )^2 ] \, du , $$
and by choosing $\epsilon$ sufficiently small we can absorb the last term on the left hand side of \eqref{eq:enestPsi1} (by taking the supremum of $t$ on the left hand side as well).

For terms $I$ and $III$ we have by \eqref{eq:deformation} that:
$$ |I| + | III | \leq \int_{\mathcal{R}_{T,U_0}} ( | L \mu | +  | X\mu | ) [ ( L \Psi_1 )^2 + ( L \Psi_2 )^2 ] \, dt du . $$
Under the bootstrap assumptions by using equation \eqref{eq:mu} we get that
\begin{equation}\label{est:muboot}
| L \mu | \leq \underline{c} \epsilon^2 , \quad | X \mu | \leq \underline{c}' \epsilon^2 , \quad \mu \geq \underline{c}'' > 0 , 
\end{equation}
for all $(t,u) \in [ 0 ,T ) \times ( - \infty , U_0 ]$, and using the above and that $T \leq \min \left( \frac{2}{\delta_1} , \frac{\min_u q^2(0,u) - q_{cr}^2}{2\delta_2} \right)$ and that neither $\delta_1$ nor $\delta_2$ are taken to be small (while $\epsilon$ is) as before, we get that
$$  |I| + | III | \leq \underline{c}'''\epsilon^2  \sup_{ t \in [ 0 ,T ]} \int_{\Sigma_t \cap \{ u \leq U_0\}}    (1+2\mu ) \mu [ (L \Psi_1 )^2 + (L \Psi_2 )^2 ] \, du , $$
and using again the smallness of $\epsilon$ we can absorb the last term on the left hand side of \eqref{eq:enestPsi1} (by taking the supremum of $t$ on the left hand side as well). So in the end it remains to show the estimates \eqref{est:muboot}.

\section{Construction of the solution in the interior region}

We work in the region 
$$\mathrm{Int} \doteq  [ 0 , T^* ) \times [ U_0 , U ] . $$
Using the bootstrap assumption \eqref{boot:non1} and using the second and fourth equation from \eqref{eq:Phi_null} we get that:
$$ | \check{\underline{L}} L \Phi_1 | (t,u) \leq \sum_{i=1}^2 \mathcal{C}' K \delta  | L \Phi_i | (t,u) , $$
and
$$   | \check{\underline{L}} L \Phi_2 | (t,u)  \leq \sum_{i=1}^2 \mathcal{C}' K \delta  | L \Phi_i | (t,u) , $$
for $(t,u) \in \mathrm{Int}$. By a matrix version of the Gr\"{o}nwall inequality and using the finiteness of the $u$-interval we get that: 
$$ | L \Phi_1 | (t,u) + | L \Phi_2 | (t,u) \leq \mathcal{C} ( U_0 , U ) ( | L \Phi_1 | ( t , U_0 ) + | L \Phi_2 | ( t , U_0 ) , $$
for $\mathcal{C} ( U_0 , U ) < \infty $ a quantity depending on $U_0$ and $U$ (that can be considered fixed as both $U_0$ and $U$ are). Using the estimates of the previous section we then get that:
$$ | L \Phi_1 | (t,u) + | L \Phi_2 | (t,u) \leq \mathcal{C}' ( U_0 , U ) \epsilon , $$
for all $(t,u) \in \mathrm{Int}$. Integrating in $t$ the last inequality we have that:
$$ | \Phi_1 | ( t,u) \leq |\Phi_1 | (0 , u ) +  \mathcal{C}' ( U_0 , U ) \epsilon \cdot t , $$
and
$$ | \Phi_2 | ( t,u) \leq |\Phi_2 | (0 , u ) +  \mathcal{C}' ( U_0 , U ) \epsilon \cdot t . $$
Note that as $t \leq \min \left( \frac{2}{\delta_1} , \frac{\min_u q^2(0,u) - q_{cr}^2}{\delta_2} \right)$ and that neither $\delta_1$ nor $\delta_2$ are taken to be small, we have for some constants $\bar{c}_1$, $\bar{c}_2$, that
$$ | \Phi_1 | ( t,u) \leq \bar{c}_1 \epsilon , $$
and
$$ | \Phi_2 | ( t,u) \leq \bar{c}_2 K , $$
for all $(t,u) \in \mathrm{Int}$. 

For the $\underline{\check{L}}$ derivatives, we use the bootstrap assumption \eqref{boot:non3} and we have that
$$ | L \underline{\check{L}} \Phi_1 | (t,u) + | L \underline{\check{L}} \Phi_2 | (t,u)  \leq 2\mathcal{C}''' K \delta \epsilon , $$
for all $(t,u) \in [0 , T^* ) \times ( - \infty , U ]$. Integrating in $t$ the last inequality we get that
$$| \underline{\check{L}} \Phi_1 | (t,u) \leq |\underline{\check{L}} \Phi_1| (0,u) + \mathrm{c} \epsilon ,$$
and
$$| \underline{\check{L}} \Phi_2 | (t,u) \leq |\underline{\check{L}} \Phi_2 | (0,u) + \mathrm{c} \epsilon ,$$
for some $\mathrm{c} > 0$ that does not depend on $K$ or $\delta$ (using that $t \leq \min \left( \frac{2}{\delta_1} , \frac{\min_u q^2(0,u) - q_{cr}^2}{2\delta_2} \right)$ and the relationship between $\delta_1$ and $K$ and $\delta$). Since $\epsilon$ is taken to be small compared to $\delta$, we get in the end that
$$ | \underline{\check{L}} \Phi_1 | (t,u) \leq \mathrm{c}_1 \delta , $$
and
$$ | \underline{\check{L}} \Phi_1 | (t,u) \leq \mathrm{c}_1 \delta , $$
for all $(t,u) \in [0 , T^* ) \times ( - \infty , U ]$. Note that we can improve the estimates above in the region $(t,u) \in [0 , T^* ) \times ( - \infty , U_c ]$ to the following:
$$ | \underline{\check{L}} \Phi_1 | (t,u) \leq \mathrm{c}_2 \epsilon , $$
and
$$ | \underline{\check{L}} \Phi_1 | (t,u) \leq \mathrm{c}_2 \epsilon . $$

\subsubsection{Formation of a shock}
Let us examine the equation for $\mu$, that is \eqref{eq:mu}. We will use the already established bounds on $\Phi_1$, $\Phi_2$, $L \Phi_1$, $L \Phi_2$, $\underline{\check{L}} \Phi_1$ and $\underline{\check{L}} \Phi_2$ to estimate the time $T^1$ where $\mu$ can potentially vanish. We will consider the case of $u \in [ U_0 , U]$ as we showed earlier that $\mu$ remains positive and bounded away from zero for $u \in (-\infty , U_0 ]$. 

Note that the estimates of the previous section can be stated as equalities in the following way:
\begin{equation}\label{eq:aux1}
 \Phi_{\kappa} (t,u ) = \Phi (0,u) + \mathcal{O} ( \epsilon ) , 
 \end{equation}
\begin{equation}\label{eq:aux2}
 \underline{\check{L}} \Phi_{\kappa} ( t,u) = \underline{\check{L}} \Phi_{\kappa} (0,u) + \mathcal{O} (\epsilon )  , 
\end{equation}
and
\begin{equation}\label{eq:aux3}
 L \Phi_{\kappa} ( t,u) =L \Phi_{\kappa} (0,u) + \mathcal{O} (\epsilon )  , 
 \end{equation}
for $\kappa \in \{1,2\}$, and for $(t,u) \in [0,T^* ) \times [U_0 , U ]$. Multiplying these two identities we get the following using the estimates we've already established:
$$ \sum_{\kappa_1 , \kappa_2 \{ 1,2\}} \Phi_{\kappa_1} (t,u) \underline{\check{L}} \Phi_{\kappa_2} (t,u) = 2 \sum_{\kappa_1 , \kappa_2 \{ 1,2\}} \Phi_{\kappa_1} (0,u) \underline{\check{L}} \Phi_{\kappa_2} (0,u) + \mathcal{O} ( K \cdot \epsilon ) + \mathcal{O} ( \delta \cdot \epsilon ) + \mathcal{O} ( \epsilon^2 ) , $$
again for $(t,u) \in [0,T^* ) \times [U_0 , U ]$. We note now that the first term of the last estimate is of size $\sim K \delta$, so all the $\mathcal{O}$ terms can be considered as errors. 

We note that 
$$ \omega_1 \sim \sum_{\kappa_1 , \kappa_2 \{ 1,2\}} \Phi_{\kappa_1}  \underline{\check{L}} \Phi_{\kappa_2} \mbox{  and  } \omega_2 \sim  \sum_{\kappa_1 , \kappa_2 \{ 1,2\}} \Phi_{\kappa_1}  L \Phi_{\kappa_2} , $$
so using equation \eqref{eq:mu}, we have that:
$$ L \mu (t,u) = \omega_1 (0,u) + \mathcal{O} ( K \cdot \epsilon ) + \mathcal{O} ( \delta \cdot \epsilon ) + \mathcal{O} ( \epsilon^2 ) . $$
We integrate this last equation in $t$ and we get that:
$$ \mu (t,u) = \mu (0,u) + \omega_1 (0,u) \cdot t + \mathcal{O} ( \epsilon ) , $$ 
using for the last terms that $t \leq   \left( \frac{2}{\delta_1} , \frac{\min_u q^2(0,u) - q_{cr}^2}{2\delta_2} \right)$. 

We take the infimum in $u$ in the last equality and we have that:
$$ \inf_{u \in [U_0 , U]} \mu (t,u) = \inf_{u \in [U_0 , U]} \mu (0,u) - \delta_1 \cdot t + \mathcal{O} (\epsilon ) . $$
Since we assume that $ \inf_{u \in [U_0 , U]} \mu (0,u) = 1$, for $\frac{2}{\delta_1} < \frac{\min_u q^2(0,u) - q_{cr}^2}{2\delta_2}$, we can verify estimate \eqref{vanishing}.

\subsubsection{Formation of the sonic line}
We will use again the established bounds on $\Phi_1$, $\Phi_2$, $L \Phi_1$, $L \Phi_2$, $\underline{\check{L}} \Phi_1$ and $\underline{\check{L}} \Phi_2$ to estimate the time $T^2$ where we can potentially have that 
$$\frac{q}{c} = 1 . $$
Recall that 
$$ \frac{q}{c} > 1 \mbox{  for $q > q_{cr}$ and } \frac{q}{c} < 1 \mbox{  for $q < q_{cr}$. } $$
We use the basic calculation:
$$ L q^2 = 2 \Phi_1 \cdot L \Phi_1 + 2 \Phi_2 \cdot L \Phi_2 , $$
and we integrate it along the integral curves of $L$, which gives us that:
$$ q^2 (t,u) = q^2(0 , u ) + 2 \int_0^t [ \Phi_1 (t',u) \cdot L \Phi_1 (t' , u) + \Phi_2 (t' , u ) \cdot L \Phi_2 (t' , u ) ] \, dt' . $$
Using equations \eqref{eq:aux1} and \eqref{eq:aux3} repeatedly, we get that
$$  q^2 (t,u) = q^2(0 , u ) + 4 (  \Phi_1 (0,u) \cdot L \Phi_1 (0 , u) + \Phi_2 (0 , u ) \cdot L \Phi_2 (0 , u ) ) \cdot t + \mathcal{O} ( \epsilon^2 ) , $$
where we used also that $t \leq   \left( \frac{2}{\delta_1} , \frac{\min_u q^2(0,u) - q_{cr}^2}{2\delta_2} \right)$. 

Using that 
$$ \Phi_2 \cdot L \Phi_2 = \mathcal{O} ( K \cdot \epsilon ) \gg \Phi_1 \cdot L \Phi_1 = \mathcal{O} ( \epsilon^2 ) , $$
from the previously established bounds on $\Phi_1$, $\Phi_2$, $L \Phi_1$ and $L \Phi_2$, we have that after taking the infimum of $q^2 (t,u)$ in $u$ the following emerges:
$$ \inf_{u \in [ U_0 , U ]} q^2 (t,u) = \inf_{u \in [ U_0 , U ]} q^2 (0,u) - 4 \delta_2 \cdot t + \mathcal{O} ( \epsilon^2 ) . $$
Setting $ \inf_{u \in [ U_0 , U ]} q^2 (t,u) = q_{cr}^2$, for $\frac{2}{\delta_1} > \frac{\min_u q^2(0,u) - q_{cr}^2}{2\delta_2}$, we can verify estimate \eqref{vanishing}.

\section{The development problem: discussion}
At this point, an interesting situation arises: according to the choice of the data, a shock might form first or a sonic line. In the later case, we cannot expect to be able to extend the solution past such a boundary, but in the former, one can try extending the solution past a shock surface. This is known as the development problem and has a long history in the area of fluid mechanics (see \cite{christ-shockdev} for an impressive result in the case of the time-dependent Euler equations, and also \cite{buckmaster-drivas-shkoller-vicol} and references therein).

In the case of hyperbolic systems in $\mathbb{R}^+ \times \mathbb{R}$ (i.e. in $1\times 1$ dimensions), there are several results on both the formation and development of shocks. Let us mention here the work of Lebaud \cite{lebaud} who showed shock formation and development for hyperbolic $p$-systems on $\mathbb{R}^+ \times \mathbb{R}$, the follow-up work of Chen and Dong \cite{chen-dong} (for the same type of equations, showing the same result for more generic data), the work of Yin \cite{yin} on shock formation and development for the compressible Euler equations with spherically symmetric data, and finally and more relevantly to our case the work of Kong \cite{kong} on shock formation and development for general quasilinear $2\times2$ hyperbolic systems on $\mathbb{R}^+ \times \mathbb{R}$. Finally, shock formation for quasilinear $N\times N$ strictly hyperbolic and genuinely nonlinear systems (see below for a definition) on $\mathbb{R}^+ \times \mathbb{R}$ has been shown in \cite{christodoulou-perez} with techniques similar to the ones used in the present paper.

The results of \cite{lebaud}, \cite{chen-dong} and \cite{kong} rely on the fact that a quasilinear $2\times 2$ system on $\mathbb{R}^+ \times \mathbb{R}$ of the form
\begin{equation}\label{eq:system2x2}
\partial_t u + A(u) \partial_x u = 0 ,
\end{equation}
for $u (t,x) = (u_1 (t,x) , u_2 (t,x) )$, where $t \in \mathbb{R}^+$, $x \in \mathbb{R}$, can be transformed into
\begin{equation}\label{eq:system2x2rm}
\begin{cases}
\partial_t v_1 + \lambda_1 ( v_1 , v_2 ) \partial_x v_1 = 0 , \\
\partial_t v_2 + \lambda_2 ( v_1 , v_2 ) \partial_x v_2 = 0 ,
\end{cases}
\end{equation}
for $\lambda_i$, $i \in \{1,2\}$, the eigenvalues of the matrix $A$. The functions $v_1$ and $v_2$ are called the \textit{Riemann invariants}.

The main result of \cite{kong} is the following:
\begin{theorem}\label{thm:kong}
Consider the system \eqref{eq:system2x2rm} in the case that it is \textit{strictly hyperbolic}, that means we assume that 
$$ \lambda_1 < \lambda_2 , $$
and \textit{genuinely nonlinear}, that means we assume that 
$$ \frac{\partial \lambda_1}{\partial v_1 } \neq 0 , \quad \frac{\partial \lambda_2}{\partial v_2 } \neq 0 . $$

Then for data 
$$ v_1 (0, x) \doteq v_{1;0} (x) , \quad v_2 ( 0, x ) \doteq v_{2;0} (x ) , \quad x \in \mathbb{R} , $$
for which there exists some $x_0 \in \mathbb{R}$ where
$$ \mbox{either  } v_{1;0}' (x_0 ) < 0 \mbox{  or  } v_{2;0}' (x_0 ) < 0 , $$
the system \eqref{eq:system2x2rm} has a unique $C^1$ solution $(v_1 , v_2 )$ in $\{ t \in [0, t_0 ] , x \in \mathbb{R} \}$ and 
$$ \limsup_{t \rightarrow t_0} \left( \| \partial_x v_1  \|_{C^0} + \| \partial_x v_2 \|_{C^0} \right) = \infty . $$

Under some additional assumptions on $\lambda_1$ and $\lambda_2$, it can be shown that the solution $(v_1 , v_2 )$ blows up at $(t_0 , x_0 )$ for $x_0 < \infty$, and moreover that in a neighbourhood $\mathcal{N}_0$ of $(t_0 , x_0 )$ we have that 
\begin{equation*}
\begin{split}
& | v_1 ( t,x) - v_1 ( t_0 , x_0 ) | \leq c_1 ( |t-t_0 | + | x - x_0 | ) \mbox{  for all $(t,x) \in \mathcal{N}_0$,} 
\\ & | \partial_x v_1 ( t,x) - \partial_x v_1 ( t_0 , x_0 ) | \leq c_1 ( |t-t_0 | + | x - x_0 | ) \mbox{  for all $(t,x) \in \mathcal{N}_0$,}
\\ & | v_2 ( t,x) - v_2 ( t_0 , x_0 ) | \leq c_2 ( |t-t_0 |^3 + | x - x_0 (t) |^2 )^{1/6} \mbox{  for all $(t,x) \in \mathcal{N}_0$,}
\\ & | \partial_x v_2 ( t,x) |  \leq c_2 ( |t-t_0 |^3 + | x - x_0 (t) |^2 )^{-1/3} \mbox{  for all $(t,x) \in \mathcal{N}_0 \setminus (t_0 , x_0 )$,}
\\ & | \partial_{xx}^2 v_2 (t,x) | \leq c_3 ( |t-t_0 |^3 + | x - x_0 (t) |^2 )^{-5/6} \mbox{  for all $(t,x) \in \mathcal{N}_0 \setminus (t_0 , x_0 )$,}
\end{split}
\end{equation*} 
for $c_1 , c_2 , c_3 > 0$ and for $x_0 (t)$ the characteristic corresponding to $\lambda_1$. 
\end{theorem}

The analysis in \cite{lebaud}, \cite{chen-dong} and \cite{kong} relies on a careful examination of the characteristics of the system \eqref{eq:system2x2rm}, that results in the treatment of a system of ODEs (something similar in spirit takes place in \cite{yin} where instead of a system provided by the Riemann invariants, an approximate system is examined following \cite{alinhac-blow1}, \cite{alinhac-blow2}). 

One can try treating the isentropic and irrotational system of compressible Euler equations in the same way. We can rewrite the equation \eqref{eq:phi} for the potential flow as a quasilinear $2\times 2$ system
\begin{equation}\label{eq:euler2x2}
\partial_y u + A_e (u) \partial_x u = 0 , 
\end{equation}
where
\begin{equation*}
A_e (u) = \begin{bmatrix} -\frac{2u_1 u_2}{c^2 - u_1^2} & \frac{c^2 - u_2^2}{c^2 - u_1^2} \\ -1 & 0 \end{bmatrix} ,
\end{equation*}
for 
$$ u = (u_1 . u_2 ) = ( \partial_y \phi , \partial_x \phi ) , $$
in the domain $\mathcal{D}$ as defined before 
$$ \mathcal{D} \doteq \{ (x, y ) \quad | \quad y \geq 0 \mbox{  for $|x| \geq 1$, and  } y \geq f(x) \mbox{  for $|x| < 1$} \}, $$
where $f$ is a smooth concave function such that 
$$ f(x) > 0 \mbox{  for $|x| < 1$,  } f(1) = f(-1) = 0 , $$
and which defines the airfoil $\mathcal{P}$ by
$$ \mathcal{P} \doteq \{ (x,y) \quad | \quad |x| < 1 , y = f(x) \} . $$
In this formulation one can set data as in Theorem \ref{thm:main}, compute the Riemann invariants (which do not seem to have a very clean form, and in the most abstract setting depend on $c$) and try to prove a result similar to \ref{thm:kong}. It should be noted though as the eigenvalues of $A_e (u)$ are given by
$$ \lambda_{\pm} (u) = \dfrac{-u_1 u_2 \pm \sqrt{ c^2 ( -c^2 + u_1^2 + u_2^2 )}}{c^2 - u_1^2} , $$
there is also the scenario of degeneration of hyperbolicity (and formation of a sonic line as in Theorem \ref{thm:main}). To this end, we are not aware of any results for quasilinear $2\times 2$ systems (only for wave equations as in \cite{jared-hyperbolic} and the references therein).

Moreover, a more interesting problem, is to try and construct (using either the first order $2\times 2$ system formulation, or the quasilinear wave equation for the potential flow) the full shock surface for appropriate data, and attempt to continue the solution past this surface for long enough time till a sonic line forms. In this direction see the notes \cite{abbrescia-speck-notes} and the references therein.


\begin{thebibliography}{10}

\bibitem{abbrescia-speck-notes}
L.~Abbrescia and J.~Speck.
\newblock The relativistic {E}uler equations: {E}{S}{I} notes on their
  geo-analytic structures and implications for shocks in $1d$ and
  multi-dimensions.
\newblock {\em ar{X}iv:2308.07289v1}, 2023.

\bibitem{alinhac-blow2}
Serge Alinhac.
\newblock Blowup of small data solutions for a class of quasilinear wave
  equations in two space dimensions {II}.
\newblock {\em Acta Math.}, 182(1):1--23, 1999.

\bibitem{alinhac-blow1}
Serge Alinhac.
\newblock Blowup of small data solutions for a quasilinear wave equation in two
  space dimensions.
\newblock {\em Ann. Math.}, 149(2):97--127, 1999.

\bibitem{transonic1}
Y.~Angelopoulos.
\newblock Instability results for transonic flows past airfoils.
\newblock {\em online preprint at ar{X}iv 2109.13369}, 2021.

\bibitem{bers-subsonic2}
Lipman Bers.
\newblock Boundary value problems for minimal surfaces with singularities at
  infinity.
\newblock {\em Transactions of the AMS}, 70:465--491, 1951.

\bibitem{bers-subsonic1}
Lipman Bers.
\newblock Existence and uniqueness of a subsonic flow past a given profile.
\newblock {\em CPAM}, 7:441--504, 1954.

\bibitem{bers-book}
Lipman Bers.
\newblock {\em Mathematical aspects of subsonic and transonic gasdynamics}.
\newblock Wiley, 1958.

\bibitem{buckmaster-drivas-shkoller-vicol}
T.~Buckmaster, T.D. Drivas, S.~Shkoller, and V.~Vicol.
\newblock Simultaneous development of shocks and cusps for 2d {E}uler with
  azimuthal symmetry from smooth data.
\newblock {\em Ann. of PDE}, 2022.

\bibitem{chen-dong}
S.~Chen and L.~Dong.
\newblock Formation and construction of shock for p-system.
\newblock {\em Science in China, Series A}, 44(9):1139--1147, 2001.

\bibitem{DC07}
D.~Christodoulou.
\newblock {\em The formation of shocks in $3$-dimensional fluids}.
\newblock Z\"urich: European Mathematical Society Publishing House., 2007.

\bibitem{christodoulou-perez}
D.~Christodoulou and D.~Raoul Perez.
\newblock On the formation of shocks of electromagnetic plane waves in
  non-linear crystals.
\newblock {\em J. Math. Phys.}, 57(081506), 2016.

\bibitem{christ-shockdev}
Demetrios Christodoulou.
\newblock {\em The shock development problem}.
\newblock Monographs in Mathematics. European Mathematical Society, 2019.

\bibitem{christ-lisibach}
Demetrios Christodoulou and Andre Lisibach.
\newblock Shock development in spherical symmetry.
\newblock {\em Ann. PDE}, 2(1):1--246, 2016.

\bibitem{christ-miao}
Demetrios Christodoulou and Shuang Miao.
\newblock {\em Compressible flow and Euler's equations}, volume~9 of {\em
  Surveys of Modern Mathematics}.
\newblock Higher Education Press, Beijing, International Press, 2014.

\bibitem{courant-friedrichs}
R.~Courant and K.O. Friedrichs.
\newblock {\em Supersonic flow and shock waves}.
\newblock Springer, 1976.

\bibitem{hksw}
G.~Holzegel, S.~Klainerman, J.~Speck, and Willie Wong.
\newblock Shock formation in small-data solutions to $3d$ quasilinear wave
  equations: An overview.
\newblock {\em J. Hyperbolic Differential Equations}, 13(1):1--105, 2016.

\bibitem{john-book}
Fritz John.
\newblock {\em Nonlinear Wave Equations, Formation of Singularities}.
\newblock University Lecture Series. AMS, 1990.

\bibitem{kong}
D.~Kong.
\newblock Formation and propagation of singularities for $2\times2$ quasilinear
  hyperbolic systems.
\newblock {\em Transactions of the AMS}, 354(8):3155--3179, 2002.

\bibitem{lebaud}
M.P. Lebaud.
\newblock Description de la formation d'un choc dans le p-syst\`{e}me.
\newblock {\em J. Math. Pures Appl.}, 73:523--565, 1994.

\bibitem{lerner-morimoto-xu}
N.~Lerner, Y.~Morimoto, and C.-J. Xu.
\newblock Instability of the {C}auchy-{K}ovalevskaya solution for a class of
  nonlinear systems.
\newblock {\em American Journal of Mathematics}, 132(1):99--123, 2010.

\bibitem{lerner-nguyen-texier}
Nicolas Lerner, Toan Nguyen, and Benjamin Texier.
\newblock The onset of instability in first order systems.
\newblock {\em J. Eur. Math. Society}, 20:1303--1373, 2018.

\bibitem{jaredjonathan}
Jonathan Luk and Jared Speck.
\newblock {Shock formation in solutions to the 2D compressible Euler equations
  in the presence of non-zero vorticity}.
\newblock {\em Invent. Math.}, 214(1):1--169, 2018.

\bibitem{majda-1}
Andrew Majda.
\newblock {\em The existence of multidimensional shock fronts}.
\newblock Mem. Am. Math. Soc. AMS, 1983.

\bibitem{majda-2}
Andrew Majda.
\newblock {\em The stability of multidimensional shock fronts}.
\newblock Mem. Am. Math. Soc. AMS, 1983.

\bibitem{metivier-nonlinear}
Guy M\'{e}tivier.
\newblock {\em Geometric Analysis of PDE and Several Complex Variables,
  Contemp. Math.}, chapter Remarks on the well-posedness of the nonlinear
  Cauchy problem, pages 337--356.
\newblock AMS, 2005.

\bibitem{morawetz1}
Cathleen Morawetz.
\newblock On the non-existence of continuous transonic flows past profiles {I}.
\newblock {\em CPAM}, IX:45--68, 1956.

\bibitem{morawetz2}
Cathleen Morawetz.
\newblock On the non-existence of continuous transonic flows past profiles
  {I}{I}.
\newblock {\em CPAM}, X:107--131, 1957.

\bibitem{morawetz3}
Cathleen Morawetz.
\newblock On the non-existence of continuous transonic flows past profiles
  {I}{I}{I}.
\newblock {\em CPAM}, XI:120--144, 1958.

\bibitem{morawetz-notes}
Cathleen Morawetz.
\newblock {\em Nonlinear waves and shocks}.
\newblock Tata Institute of Fundamental Research, 1981.

\bibitem{morawetz-ams}
Cathleen Morawetz.
\newblock The mathematical approach to the sonic barrier.
\newblock {\em Bulletin of the AMS}, 6(2):127--145, 1982.

\bibitem{shiffman}
Max Shiffman.
\newblock On the existence of subsonic flows of a compressible fluid.
\newblock {\em Journal of Rational Mechanics and Analysis}, 1:605--652, 1952.

\bibitem{speck-notes}
Jared Speck.
\newblock Lecture notes on shock formation in quasilinear wave equations: {A}n
  overview of the nearly plane symmetric gime.
\newblock {\em Unpublished notes}.

\bibitem{speck-book}
Jared Speck.
\newblock {\em Shock Formation in Small-data Solutions to 3{D} Quasilinear Wave
  Equations}.
\newblock American Mathematical Society, 2016.

\bibitem{jared-hyperbolic}
Jared Speck.
\newblock Finite-time degeneration of hyperbolicity without blowup for
  quasilinear wave equations.
\newblock {\em Ann. PDE}, 10(8):2001--2030, 2018.

\bibitem{jaredgustavjonwillie}
Jared Speck, Gustav Holzegel, Jonathan Luk, and Willie Wong.
\newblock Stable shock formation for nearly simple outgoing plane symmetric
  waves.
\newblock {\em Ann. PDE}, 2(2):1--198, 2016.

\bibitem{yin}
H.~Yin.
\newblock Formation and construction of a shock wave for 3-{D} compressible
  {E}uler equations with the spherical initial data.
\newblock {\em Nagoya Math J.}, 175:125--164, 2004.

\end{thebibliography}
\end{document}